\input amssym.def
\input amssym.tex

\font\twelverm=cmr12 at 14pt
\hsize = 31pc
\vsize = 45pc
\overfullrule = 0pt

\def\www{\hat}
\def\bd{\buildrel .\over +}
\def\vs{\vskip 1pc}
\font\ninerm=cmr10 at 9pt

 \font\newrm =cmr10 at 24pt
\def\bul{\raise .9pt\hbox{\newrm .\kern-.105em } }

 \def\fr{\frak}

\baselineskip=13pt
 
 \def\h{\hbox{ }}

 \def\u{{\fr u}}

 \def\n{{\fr n}}
 \def\a{{\fr a}}
 \def\d{{\fr d}}

 \def\k{{\fr k}}
 \def\b{{\fr b}}
 
 \def\hh{{\fr h}}
 \def\tt{{\fr t}}

 \def\I{{\cal I}}
 \def\g{{\fr g}}
 \def\v{{\fr v}}

 \def\<{\le}
 \def\>{\ge}

 \def\s{{\h\subset\h}}
 
 \def\vs{\vskip }

 \def\mapright#1
  {\smash{\mathop
  {\longrightarrow}
  \limits^{#1}}}

 \def\kk#1{{\kern .4 em} #1}
 \def\vs{\vskip 1pc}

\font\smallbf=cmbx10  at 9pt 

\font\authorfont=cmr12 at 12pt
\font\ninerm=cmr10 at 9pt
\centerline{\twelverm Powers of the Euler product and commutative subalgebras}
\vskip 2pt
 \centerline {\twelverm of a complex simple Lie algebra} 
\vskip 1.5pc
\baselineskip=11pt
\vskip8pt
\centerline{\authorfont BERTRAM KOSTANT\footnote*{\ninerm
Research supported in part by NSF contract DMS-0209473 and the KG\&G Foundation.}}\vskip 2pc

 \baselineskip=12pt 
\noindent{\smallbf ABSTRACT.} \ninerm If $\frak g$ is a complex simple Lie algebra, and $k$ does not exceed the dual Coxeter
number of
$\frak g$, then the k$^{th}$ coefficient of the $dim\,\, \frak g$
power of the Euler product may be given by the dimension of a subspace of $\wedge^k\frak g$ defined by all
abelian subalgebras of $\frak g$ of dimension $k$. This has implications for all the coefficients of all
the powers of the Euler product. Involved in the main results are Dale Peterson's $2^{rank}$ theorem
on the number of abelian ideals in a Borel subalgebra of $\frak g$, an element of type $\rho$ and my
heat kernel formulation of Macdonald's $\eta$-function theorem, a set
$D_{alcove}$ of special highest weights parameterized by all the alcoves in a Weyl chamber
(generalizing Young diagrams of null $m$-core when $\frak g= Lie\,Sl(m,\Bbb C)$), and the homology and
cohomology of the nil radical of the standard maximal parabolic subalgebra of the affine Kac-Moody
Lie algebra.

\vskip 1pc 

\centerline{\bf 0. Introduction}\vskip 1pc

\rm
 0.1. The Euler product
is the power series in the variable $x$ obtained from the expansion
of the infinite product
$\Pi_{n=1}^{\infty}(1-x^n)$. The Dedekind $\eta$-function is
$x^{1/24}$ times the Euler product. Let
$\g$ be a complex simple Lie algebra, and let $K$ be a simply connected compact group such that 
$\k = Lie\,K$ is a compact form of $\g$. Let $\ell$ be the rank of $K$ and let $T\s K$
be a maximal torus of $K$. Macdonald has given a
formula for the $dim\,K$ power of the Euler product in terms of a
summation over a lattice in $\hh$. See [Ma-1]. A new understanding
of this formula was made by V. Kac. It arose from his dominator identity. See Chapter
12 in [Ka]. Another approach to the formula was a consequence of
the Laplacian determination (an infinite dimensional analogue of Theorem 5.14
in [Ko-3]) of the homology of the 
``nilradical" of the standard maximal parabolic
subalgebra of the affine Kac-Moody Lie algebra. This is due to H. Garland in [G]. Because the 
Laplacian is positive semidefinite this
approach has the advantage of implying some very important inequalities. Garland's result is a
special case of a much more general result of Shrawan Kumar. See [Ku-1] (or
Theorem 3.4.2 in [Ku-2]). Kumar's Theorem 3.5.2 is a far reaching infinite
dimensional analogue of Theorem 5.7 in [Ko-1].

Let $\hh = i\,Lie\,T$ and identify $\hh$ with its dual using the Killing form so that
$\Delta\s \hh$ where $\Delta$ is the set of roots for the pair
$(\hh_{\Bbb C},\g)$. Let $\hh^+\s \hh$ be a Weyl chamber defined by
a choice, $\Delta_+$, of positive roots and let $D\s \hh^+$ be the
set of dominant integral forms on $\hh$. For each $\lambda\in D$ let
$\pi_{\lambda}:K\to Aut\,V_{\lambda}$ be an irreducible
representation with highest weight $\lambda$. Let $\chi_{\lambda}$ be
the character of $\pi_{\lambda}$ and let $Cas(\lambda)$ be the scalar
value taken by the Casimir element (relative to the Killing form)
on $V_{\lambda}$. Let $\rho$ be one half the sum of the positive
roots and let $a_P = exp\,2\,\pi\,i\,\,2\rho$. An element in $K$
conjugate to $a_P$ is referred to in [Ko-3] as an element of type
$\rho$. Using Macdonald's formula the
following result was established as part of Theorem 3.1 in [Ko-3].
\vs {\bf Theorem 0.1.} {\it For any
$\lambda\in D$ one has $\chi_{\lambda}(a_P)\in \{-1,0,1\}$ and 
$$(\prod_{n=1}^{\infty} (1-x^n))^{dim\,K} =
\sum_{\lambda\in
D}\chi_{\lambda}(a_P)\,\,dim\,V_{\lambda}\,\,x^{Cas(\lambda)}\eqno
(0.1)$$} \vs For a heat kernel (on $K$) interpretation of (0.1) see
\S 5 in [Ko-3] or [F]. \vs Of course the only dominant weights $\lambda$
which contribute to the sum (0.1) are those in the subset $\{\lambda\in D\mid
\chi_{\lambda}(a_P)\in \{-1,1\}$. The determination of this subset is implicit
in Lemma 3.5.2 of [Ko-3]. This is clarified in \S 2.2 of the present paper. In more
detail let $W_f$ be the affine Weyl group operating in $\hh$ and let $A_1=\{x\in
\hh^+\mid \psi(x)\leq 1\}$ where $\psi\in \Delta_+$ is the highest root. Then $A_1$ is
a fundamental domain for the action of $W_f$ and any subset of form $A_{\sigma} =
\sigma(A_1)$, for $\sigma\in W_f$, is referred to as an alcove. Let $W_f^+=
\{\sigma\in W_f\mid A_{\sigma}\s \hh^+\}$. An alcove $A_{\sigma}$ is called dominant
if $\sigma\in W_f^+$. (A study of the set
 of dominant alcoves was made by Bott in [B] in
connection with the topology of the loop group
$\Omega(K)$). The set
of dominant alcoves, or rather $W_f^+$, parameterize a subset $D_{alcove}$ of $D$ by
 defining $\lambda^{\sigma} = \sigma(2\rho)/2 - \rho$, for any $\sigma\in W_f^+$, and
putting $D_{alcove} = \{\lambda^{\sigma}\mid \sigma\in W_f^+\}$. The element $2(\rho +
\lambda)$ is in the interior of $A_{\sigma}$ and in fact, if $\g$ is simply-laced, this
element is the center of the inscribed sphere of the simplex $A_{\sigma}$. (See the
cautionary statement about $\lambda^{\sigma}$ in Remark 3.4). Let $b_k,\,k\in\Bbb Z_+$,
be the coeficients defined by the expansion $(\prod_{n=1}^{\infty} (1-x^n))^{dim\,K} =
\sum_{k=0}^{\infty} b_kx^k$. One then has (taken from Theorem 2.4 and (3.40))\vs {\bf
Theorem 0.2.} {\it Let
$\lambda\in D$. Then $\chi_{\lambda}(a_P)\in \{-1,1\}$ if and only if $\lambda\in
D_{alcove}$ so that in such a case $\lambda = \lambda^{\sigma}$ for a unique $\sigma\in
W_f^+$. Moreover in such a case $\chi_{\lambda^{\sigma}}(a_P)=(-1)^{\ell(\sigma)}$ where
$\ell(\sigma)$ is the length of $\sigma$. In addition $Cas(\lambda^{\sigma})\in \Bbb Z_+$
and one has the finite sum $$b_k = \sum_{\sigma\in W_f^+,\,\, Cas(\lambda^{\sigma}) = k}
(-1)^{\ell(\sigma)}\,dim\,V_{\lambda^{\sigma}}\eqno (0.2)$$}\vs {\bf Remark 0.3.} If $K
= SU(m)$ then as one knows the representation theory of $SU(m)$ defines a bijection
$$D\to P\eqno (0.3)$$ where $P$ is the set of partitions of length at most $m-1$.
 For any
$p\in P$ one defines another partition $\tilde p$ called its $m$-core. One
says that
$p$ has null $m$-core if $\widetilde {p}$ is the empty partition. From the first
statement in Theorem 0.2 it follows from Exercise 1.1.8(c) in [M] or \S 3.4 in [A-F]
or p. 467-469 in [St] that the image of $D_{alcove}$ in (0.3) is exactly the set of
$p\in P$ with empty $m$-core. Such a partition exists only if the size of $p$ is an
integral multiple $mk$ of $m$ and using, for example results of Bott, one can show, in
this case, that the number of such partitions is $m + k -2\choose m-2$. 
\vs 0.2. Let $\www\g$ be the affine Kac-Moody Lie algebra corresponding to $\g$. 
The ``nil radical" of a standard maximal parabolic
subalgebra is either $\u^- = t^{-1}\,\g[t^{-1}]$ or the opposed algebra $\u = t\,\g[t]$.
The exterior algebra $\wedge\u^-$ is bigraded by $\Bbb Z_+\times \Bbb Z_+$ with
homogeneous components
$(\wedge^n\u^-)_k$ where $-k$ is the $t$-degree. Furthermore $\wedge\u^-$ is an
underlying space for the chain complex whose derived homology is $H_*(\u^-)$. In
addition it is also the underlying space for a cochain complex whose derived cohomology
${\cal H}(\u)$ is a restricted form (see \S 4.3) of $H^*(\u)$. Garland's theorem
(Theorem 3.10 here), Theorem 3.13 and Theorem 4.8 yield\vs {\bf Theorem 0.4.} {\it As
$\g$-modules, $H_*(\u^-)$ and $H^*(\u)$ are equivalent and multiplicity free. In fact
$$\eqalign{H_*(\u^-)&\equiv {\cal H}(\u)\cr &\equiv \sum_{\sigma\in
W_f^+}V_{\lambda^{\sigma}}\cr}\eqno (0.4)$$ With respect to the bigrading
$$\eqalign{(H_n(\u^-))_k&\equiv ({\cal H}^n(\u))_k \cr &\equiv \sum_{\sigma\in
W_f^+,\,\ell(\sigma) = n,\,\,Cas(\lambda^{\sigma})= k}V_{\lambda^{\sigma}}\cr}\eqno
(0.5)$$}\vskip .5 pc 0.3. For $k\in \Bbb Z_+$ let $C_k\s \wedge^k\g$ be the span of all
1-dimensional subspaces of the form $\wedge^k\a$ where $\a\s \g$ is any $k$-dimensional
abelian subalgebra of $\g$. Let $M$ be the maximal dimension of a commutative
subalgebra of $\g$. Obviously $C_k \neq 0$ if and only if $k\leq M$. The value of $M$
was computed by Malcev for each $\g$-type (e.g. $M= 36$ if the $\g$ type is
$E_{8}$).

If $\v\s
\g$ is any (complex) subspace of
$\g$ which is stable under $ad\,\hh$, let $\Delta(\v) = \{\varphi\in \Delta\mid
e_{\varphi}\in
\v\}$ where
$e_{\varphi}$ is a root vector for $\varphi$. Let $\b\s \g$ be the Borel subalgebra
containing $\hh$ such that $\Delta(\b) = \Delta_+$ and let $\n= [\b,\b]$. Let $\Xi$ be an
index set parameterizing the set of all abelian ideals of $\b$ and for any $\xi\in
\Xi$ let $ \a_{\xi}$ be the corresponding ideal. It is immediate that $\a_{\xi}\s \n$.
Let $$\lambda_{\xi} = \sum_{\varphi\in \Delta(\a_{\xi})} \varphi\eqno (0.6)$$ Let $\Xi_k =
\{\xi\in \Xi\mid dim\,\a_{\xi} = k\}$. If $\xi\in \Xi_k$ let $V_{\xi}\s \wedge^k\g$ be
the $\g$-submodule generated by $\wedge^k\a_{\xi}$ with respect to the adjoint action of
$\g$ on $\wedge\g$. Obviously $V_{\xi}\s C_k$. Furthermore it is immediate that $V_{\xi}$
is irreducible,
$\wedge^k\a_{\xi}$ is the highest weight space in $V_{\xi}$ and $\lambda_{\xi}$ is the
highest weight of
$V_{\xi}$. Moreover we have proved (see Theorem 4.2) the following result as part of
Theorem (8) in [Ko-3]. \vs {\bf Theorem 0.5.} {\it For any $\k\in \Bbb Z_+$ where
$0\leq k\leq M$ one has the direct sums 
$$C_k = \sum_{\xi\in\, \Xi_k} V_{\xi}\eqno (0.7)$$ and $$C = \sum_{\xi\in
\,\Xi} V_{\xi}\eqno (0.8)$$ Furthermore $C$ and, a fortiori, $C_k$, are
multiplicity-free $\g$-modules.}\vs It is a beautiful later result of Dale Peterson that
$card\,\Xi = 2^{\ell}$. A simpler proof of Peterson's theorem was given in [C-P]. See
Theorem 2.9 in that reference. It is clear that there are $2^{\ell}$ alcoves in $2\,A_1$.
They are parameterized by $W_f^{(2)} = \{\sigma\in W_f^+\mid A_{\sigma}\s 2\,A_1\}$.
The Cellini-Papi proof of Peterson's theorem is obtained by establishing a bijection
$\Xi\to W_f^{(2)}$. The bijection is stated in a new way here (Theorem 0.6 below ) using
$D_{alcove}$ and the weights $\{\lambda_{\xi}\mid
\xi\in \Xi\}$. Recently R. Suter in [Su] has shown the Peterson's theorem follows
from Theorem (8) in [Ko-3]. With the benefit of this knowledge we have found a new proof
of Theorem 0.6 without the use of [C-P]. The following statement is Theorem 4.3. \vs
{\bf Theorem 0.6.} {\it For any
$\xi\in
\Xi$ there exists an (necessarily unique) element $\sigma_{\xi}\in W_f^+$ such that $$\lambda_{\xi} 
= \lambda^{\sigma_{\xi}}\eqno (0.9)$$ Moreover $\sigma_{\xi}\in W_f^{(2)}$ and the map $$\xi\to
W_f^{(2)},\qquad \xi\mapsto \sigma_{\xi}\eqno (0.10)$$ is a bijection. In particular one has
the inclusion $$\{\lambda_{\xi}\mid \xi\in \Xi\} \s D_{alcove}\eqno (0.11)$$ and the count
(Peterson's theorem) $$card\,\Xi = 2^{\ell}\eqno (0.12)$$}\vs {\bf Remark 0.7.} Suter in [Su] has
independently discovered (0.9).\vs The $2^{\ell}$ element subset $\{\Lambda^{\sigma}\mid \sigma\in
W_f^{(2)}\}$ of $D_{alcove}$ is characterized in (see Theorem 4.5)  \vs {\bf Theorem 0.8.} {\it
Let $\sigma\in W_f^+$. Then $$Cas(\lambda^{\sigma}) \geq \ell(\sigma)\eqno (0.13)$$ and equality
occurs in (0.13) if and only if $\sigma\in W_f^{(2)}$. Furthermore in that case writing $\sigma =
\sigma_{\xi}$ for $\xi\in \Xi$ (Theorem 0.6) one has $$\eqalign{Cas(\lambda^{\sigma}) &=
\ell(\sigma)\cr &= dim\,\a_{\xi}\cr}\eqno (0.14)$$}\vs Cup product ($\vee$) defines the structure
of an algebra on ${\cal H}(\u)$. Theorem 0.8 has implications on the determination of an important
subalgebra of ${\cal H}(\u)$. Following a suggestion of Pavel Etingof, introduce
a new grading, ${\cal H}^{[j]}(\u),\,j\in \Bbb Z_+$, in ${\cal H}(\u)$ by putting ${\cal
H}^{[j]}(\u) = \sum_{n,k\in \Bbb Z,\,k-n = j} ({\cal H}^n(\u))_k$. The homogeneous components
are finite dimensional and one has $${\cal H}(\u) = \sum_{j=0}^{\infty}{\cal H}^{[j]}(\u)$$ and
$${\cal H}^{[j]}(\u)\vee {\cal H}^{[j']}(\u) \s {\cal H}^{[j+j']}(\u)$$ In particular ${\cal
H}^{[0]}(\u)$ is a finite dimensional subalgebra of ${\cal H}(\u)$. See Proposition 4.12.
Etingof suggested that our results should yield the structure of ${\cal H}^{[0]}(\u)$.
Indeed this the case. Identify
$\g$ with its dual using the Killing form. Then
$d:\g\to
\wedge^2\g$ where $d$ is the Cartan-Eilenberg-Koszul coboundary operator whose derived cohomology
is $H^*(\g)$. Let $(d\,\g)$ be the ideal in $\wedge\g$ generated by $d\,\g$. Theorem 4.3 of
[Ko-6] establishes that one has the following direct sum $$\wedge\g = C\oplus (d\,\g)\eqno
(0.15)$$ so that $C$ inherits an algebra structure. I had no idea of the meaning of this algebra
when [Ko-6] was written. Its meaning is clarified in (see Theorem 4.16 in the present paper) \vs {\bf
Theorem 0.9.} {\it As a $\g$-module ${\cal H}^{[0]}(\u)$ is
multiplicity-free with $2^{\ell}$ irreducible components. In fact (see (0.8)) 
$${\cal H}^{[0]}(\u)
\equiv \sum_{\xi\in \Xi} V_{\xi}\eqno (0.16)$$ As an algebra (under cup product)
$${\cal H}^{[0]}(\u)\equiv \wedge\g/ (d\g)\eqno (0.17)$$} \vs Let $(d\g)^k = (d\g)\cap \wedge^k\g$.
Also let $\theta(Cas)\in End\,\wedge\g$ be the action of $Cas$ induced by the adjoint
representation of $\g$ on $\wedge\g$. Then (2.1.7) and Theorem (5) in [Ko-2] yield \vs {\bf Theorem
0.10.} {\it Let $k\in \Bbb Z_+$. Then the following four numbers are all equal
$$\eqalign{&[1]\,\,dim\,C_k\cr&[2]\,\, dim\, \wedge^k\g/(d\g)^k \cr &[3]\,\,dim\,\{v\in
\wedge^n\g|\theta(Cas)\,v = k\,v\}\cr &[4]\,\,dim\,({\cal H}^k(\u))_k\cr}$$}\vskip .5pc 0.4. One
difficulty in using (0.2) to compute the coefficient $b_k$ is the alternation in signs in (0.2).
By Theorem 0.8 this difficulty would disappear if $k$ were such that $Cas(\lambda^{\sigma}) = k$
implies that $\sigma\in W_f^{(2)}$. But this is the case if $k\leq h^{\vee}$ where $h^{\vee}$ is
the dual Coxeter number. The following is one our main results (see Theorem 4.23). \vs {\bf
Theorem 0.11.} {\it Let
$k\leq h^{\vee}$. The following seven number are all equal. $$\eqalign{&[1]\,\,(-1)^kb_k\cr 
&[2]\,\,dim\, C_k \cr
&[3]\,\, \sum_{\xi \in \Xi_k}dim\, V_{\xi} \cr
&[4]\,\,dim\, \{v \in \wedge^k \g\mid \theta(Cas)\,v = k\, v\} \cr
&[5]\,\,dim\, \wedge^k\g/(d\g)^k\cr
&[6]\,\,dim\,{\cal H}^k(\u)\cr
&[7]\,\, dim\, H_k(\u^-)\cr}$$}\vs {\bf Example}. If $K = SU(5)$ then since $dim\,SU(5) = 24$ one
has $b_n = \tau(n+1)$ where $n\mapsto \tau(n)$ is the Ramanujan tau function. Here $h^{\vee}
= 5$. Theorem 0.11 says $\wedge\,Lie\,\, Sl(5,\Bbb C)$ ``sees" the first five non-trivial Ramanujan
numbers. One has $\tau(2) = -24,\,\,\tau(3) = 252,\,\,\tau(4) = -1472,\,\,\tau(5) = 4830\,\,\tau
(6)= -6048$. See [Se], p.97. One readily checks that (choosing [4] in Theorem
0.11 for the computation), $$\eqalign{dim\,C_1&= 24\cr dim\,C_2&= 252\cr dim\,C_3&= 1472\cr dim\,C_4
&= 1472\cr dim\,C_5&= 6048\cr}$$\vskip .5pc P. Etingof points out that Theorem 0.11 can be expressed
as a homology acyclicity statement. Let $\partial_-$ be the boundary operator on $\wedge\u^-$ whose
derived homology is $H_*(\u^-)$. The restriction of $\partial_-$ to $(\wedge\u^-)_k$, for any $k\in
\Bbb Z_+$, defines a finite dimensional
 subcomplex (i.e. the $t$-degree is fixed to be $-k$)
$$(\wedge^k\u^-)_k
\longrightarrow (\wedge^{k-1}\u^-)_k\longrightarrow\cdots \longrightarrow (\wedge^{0}\u^-)_k
\longrightarrow 0\eqno (0.18)$$ We thank Etingof for the following statement (see Theorem 4.24).
\vs {\bf Theorem 0.12.} {\it If
$k\leq h^{\vee}$ then the complex (0.18) is acyclic. That is, $(H_n(\u^-))_k = 0$ unless $n=k$ so
that 
$(H_*(\u^-))_k = (H_k(\u^-))_k$. In fact $(H_*(\u^-))_k = H_k(\u^-)$ and hence $$dim\,(H_*(\u^-))_k =
(-1)^kb_k\eqno (0.19)$$}\vskip .5pc 0.5. One can raise the Euler product to the $s$ power where $s$
is any complex number, by taking its logarithm, multiplying by $s$ and exponentiating. One then
easily has that $$(\Pi_{n=1}^{\infty}(1-x^n))^s = \sum_{k=0}^{\infty}f_k(s)\,x^k\eqno (0.20)$$ where
$f_k(s)$ is a polynomial of degree $k$. Obviously $$b_k = f_k(dim\,K)\eqno (0.21)$$ Although one can
give an expression for these polynomials the expression yields very little understanding of the
nature of the polynomials. One approach could be a determination of the roots of the $f_k$. It
is a long standing question (see p. 98 in [Se]) about the Ramanujan numbers as to whether 24 is
ever a root of the $f_k$. The following result suggests a possible Lie-theoretic connection with
the roots. Obviously $f_1 =1$ and $0$ is a root of $f_k$ for all $k>0$. Consider $f_2,f_3$ and
$f_4$. Since 4 is neither a pentagonal number nor a triangular number it follows from a formula
(s=1) of Euler and (s=3) of Jacobi that for one missing root $r_4$ one has $f_4(s) = 1/4!
s(s-1)(s-3)(s-r_4)$. Similarly $r_3$ and $r_2$ exist so that $-f_3(s) = 1/3!s(s-1)(s-r_3)$ and
$f_2(s) = 1/2! s(s-r_2)$. On the other hand the only cases where $M<h^{\vee}$ are when $\g$ is of
type $A_1,A_2$ and $G_2$. As a consequence of Theorem 0.11 one has (Theorem 4.27)\vs {\bf Theorem
0.13.} {\it The missing roots $r_4,\,r_3$ and $r_2$ are,
respectively, the complex dimensions of $G_2$, $A_2$ and $A_1$, namely 14,8 and 3 so that 
$$\eqalign{f_4(s) &= 1/4!\,\,s(s-1)(s-3)(s-14)\cr -f_3(s) &= 1/3!\,\,s(s-1)(s-8)\cr f_2(s) &= 1/2!\,\,
s(s-3)\cr}$$}\vskip .5pc Let $k$ be any positive integer. If $m\in \Bbb Z_+$ and $m\geq 2$ let
$C_k(m)$ equal
$C_k$ for the case where $K= SU(m)$. If $m\geq k$ then $k\leq h^{\vee}= m$ and hence, by Theorem
0.11,
$$f_k(m^2-1) =(-1)^k dim\,C_k(m)\eqno (0.22)$$ The following result implies that $f_k(s)$ is encoded in
the $k$-dimensional commutative subalgebra structure of $Lie\,Sl(m,\Bbb C)$ for $k$ different values
of
$m$ where
$m\geq k$ and
$m>1$ (see Theorem 4.28).
\vs {\bf Theorem 0.14.} {\it Let $k$ be a positive integer. Then $f_k(s)$ is determined by the
numbers $dim\,C_k(m)$ for $k$ different values of $m\in \Bbb Z_+$ where $m\geq k$ and $m>1$.}\vskip .5
pc 0.6. We wish to thank Pavel Etingof, Shrawan Kumar and Richard Stanley for valuable
informative conversations.  

\vs

\centerline{\bf 1. Alcoves and the affine Weyl group}\vs \rm 1.1.  Let $\g$ be a complex simple
Lie algebra. The value of the Killing form
$B$ of $\g$ on
$x,y\in
\g$ will be denoted by $(x,y)$ and using $B$ we will
identify $\g$ with its dual space. Let 
$\k$ be a compact form of $\g$ and let $K$ be a
corresponding simply connected compact Lie group. Let $T$ be a maximal torus of $K$ and let $\tt=
Lie\,T$. Let $\tt_{\Bbb C}\s \g$ be the complexification of $\tt$ so that $\tt_{\Bbb C}$ is a Cartan
subalgebra of $\g$. The restriction of $B$ to
$\tt_{\Bbb C}$ will be used to identify $\tt_{\Bbb C}$ with its dual
space so that
$\Delta \s
\tt_{\Bbb C}$ where $\Delta$ is set of roots for the pair
$(\tt_{\Bbb C,}\g$). The span of $\Delta$ over $\Bbb R$ is a real form $\hh$ of $\tt_{\Bbb C}$. In
fact $\hh = i\tt$ and one knows that $B|\hh$ is positive definite. 

Let $\ell = \hbox{rank}\,\k$ and
let $I = \{1,\ldots,\ell\}$. Let
$\Delta_+\s \Delta$ be a choice of a positive root system and let
$\Pi$ be the set of simple positive roots. The elements of $\Pi$
will be indexed by $I$ so that we can write $\Pi =
\{\alpha_i\},\,i\in I$. Let $\varepsilon$ be the epimorphism
$$\varepsilon:\hh\to T\eqno (1.1)$$ where $\varepsilon(x) =
exp\,\,2\pi i\,\,x$. Let $\Gamma$ be the kernel of $\varepsilon$ 
so that
$\Gamma$ is a lattice in $\hh$. For any $\varphi\in \Delta$ one
knows that $$\varphi^{\vee}\in \Gamma\eqno (1.2)$$ where
$\varphi^{\vee} = 2\,\varphi/(\varphi,\varphi)$. Furthermore the
set of elements $\{\alpha_i^{\vee}\},\,i\in I$ is a basis of
$\Gamma$ so that any
$\gamma\in \Gamma$ can be uniquely written $$\gamma= \sum_{i\in
I}m_i\,\, 2\,\alpha_i/(\alpha_i,\alpha_i)\eqno (1.3)$$

1.2. For any
$z\in \hh$ let $t_z:\hh \to \hh$ be the translation map by $z$ so
that $t_z(x) = z +x$ for any $x\in \hh$. Let 
$W$ be the Weyl group for the pair
$(T,K)$. Obviously $\Gamma$ is stabilized by $W$ with respect
to the action of
$W$ on
$\hh$. Consequently the translation group $\widetilde {\Gamma} =
\{t_{\gamma}\mid \gamma \in \Gamma\}$ is normalized by $W$.
 The affine Weyl group $W_f$ is the semidirect product 
 $$W_f = \widetilde {\Gamma} \rtimes W$$ and we will be mainly concerned with its natural affine
action on $\hh$. 

For any
$n\in \Bbb Z$ and
$\varphi\in \Delta_+$ let $\hh_{\varphi,n}$ be the hyperplane in
$\hh$ defined by putting $$\hh_{\varphi,n}= \{x\in
\hh_{\Bbb R}\mid (\varphi,x)= n\}$$ We will use the
word wall to refer to a hyperplane in $\hh$ of this form. More specifically a wall of the form
$\hh_{\varphi,n}$ will referred to as a $\varphi$-wall. An element
$x\in
\hh$ will be called $W_f$-singular if $x\in Sing(\hh)$ where 
$$Sing(\hh) = \bigcup_{\varphi\in \Delta_+,\,n\in \Bbb
Z}\hh_{\varphi,n}$$ An element $y\in \hh$ will be called
$W_f$-regular if $y$ lies in the complement $Reg(\hh)$ of
$Sing(\hh)$ in $\hh$. The closure $A$ of a connected component of
$Reg(\hh)$ is called an alcove. The connected component itself is
clearly $Reg(A)$ where $Reg(A) = A\cap Reg(\hh)$ and one readily has that
$Reg(A)$ is the interior of $A$. The affine Weyl group
$W_f$ operates simply and transitively on the set ${\cal A}$ of all
alcoves. Let
$\psi\in \Delta_+$ be the highest root. A special alcove $A_1$,
referred to as the fundamental alcove, can be defined by $$A_1 =
\{x\in\hh\mid (\alpha_i,x)\geq 0,\,i\in
I,\,\hbox{and}\,\,(\psi,x)\leq 1\}$$ We can then index the elements
of
${\cal A}$ by $W_f$ where, if $\sigma\in W_f$, we put $A_{\sigma} =
\sigma(A_1)$. 

Let $\sigma\in W_f$ and put $T_{\sigma} =
\varepsilon(A_{\sigma})$. One knows (1) that every element in $\hh$
is $W_f$-conjugate to a unique element in $A_1$,
 (2) $\varepsilon:W_1\to T_1$ is bijective (see
(1.1)), and (3) any element in $K$ is $K$-conjugate to a unique
element in
$T_1$. Since these properties are clearly preserved by the action of
$W_f$ one immediately has \vskip 1pc {\bf Proposition 1.1.} {\it Let
$\sigma\in W_f$. Then (1) every element in $\hh$
is $W_f$-conjugate to a unique element in $A_{\sigma}$, (2)
$\varepsilon:A_{\sigma}\to T_{\sigma}$ is bijective (see (1.1)), and
(3) any  element in $K$ is $K$-conjugate to a unique
element in $T_{\sigma}$. }  \vskip 1pc  1.3. For any $(\varphi,n)
\in \Delta_+\times \Bbb Z$ let $s_{\varphi,n}$ be the reflection
in $\hh$ defined by the wall $\hh_{\varphi,n}$. We write
$s_{\varphi} = s_{\varphi,0}$. Of course $s_{\varphi} \in W$.
However for any $n\in \Bbb Z$ one readily sees that $$s_{\varphi,n}
= t_{n\,\varphi^{\vee}}\,\,s_{\varphi}\eqno (1.4)$$ so that
$s_{\varphi,n}\in W_{f}$. In fact one knows that $W_f$ is a Coxeter
group with the
$\ell + 1$ generators $\{s_i,\,s_{\psi,1}\}, i\in I$, where we have
written $s_i = s_{\alpha_i}$. In particular one has a length
function, $\sigma\mapsto \ell(\sigma)$ on $W_f$.
For any $\varphi\in \Delta_+$ and $\sigma\in W_f$ let $$n_{\varphi}(\sigma) = 
\#\,\,\hbox{of $\varphi$-walls separating $Reg(A_{\sigma})$ from $Reg(A_1)$}\eqno (1.5)$$ 
It follows easily that if $\sigma\in W_f$ then $\ell(\sigma)$ can be given
geometrically by 
$$ \eqalign{\ell(\sigma) &=
\#\,\,\hbox{of walls separating $Reg(A_{\sigma})$ from $Reg(A_1)$}\cr
&=\sum_{\varphi\in \Delta_+}n_{\varphi}(\sigma)\cr} \eqno (1.6)$$\vskip 1pc We adopt the convention
that $\Bbb Z_+$ is the set of non-negative integers and $\Bbb N$ is the set of positive integers.\vs
 {\bf
Remark 1.2}. Note that for any $\sigma\in W_f$ and any
$\varphi\in \Delta_+$ there exists $n\in \Bbb Z$ such that for the open and closed unit intervals,
$(n,n+1)$ and $[n,n+1]$ in
$\Bbb R$, one has
$$\eqalign{\varphi(Reg(A_{\sigma}))
&\s (n,n+1)\cr \varphi(A_{\sigma}) &\s
[n,n+1]\cr} \eqno (1.7) $$ This is immediate since it is clearly true for $A_{\sigma} = A_1$. One
notes also that if $n\in \Bbb Z_+$ then $$n= n_{\varphi}(\sigma)\eqno (1.8)$$
 \vs
 1.4.  Let $\hh^+\s \hh$ be 
the Weyl chamber corresponding to
$\Delta_+$ so that $$\hh^+ = \{x\in \hh\mid (\varphi,x)\geq
0,\,\forall
\varphi\in
\Delta_+\}$$ The interior $Int(\hh^+)$ can be characterized by
$$Int(\hh^+) = \{x\in \hh\mid (\alpha_i,x)>0,\,\forall i\in I\}\eqno
(1.9)$$ Let ${\cal A}^+$ be set of all alcoves $A$ such that $A\s
\hh^+$. This defines a subset $W_f^+$ of the affine Weyl group by the condition
${\cal A}^+= \{A_{\sigma}\mid \sigma\in W_f^+\}$. Note that by
Remark 1.2 and (1.9) one has $\sigma\in W_f^+$ if and only if
$Int(\hh^+)\cap A_{\sigma} \neq \emptyset$. It follows easily then
that
$$\hh^+ = \bigcup_{\sigma\in W_f ^+} A_{\sigma}\eqno (1.10)$$\vskip
.5pc {\bf Remark 1.3.} One readily shows that $W_f^+$ is the set of
minimal length representatives of the 
right cosets of $W$ in $W_f$. In fact if $w\in W$ and $\sigma\in W_f^+$ then $$\ell(w\,\sigma) =
\ell(w) + \ell(\sigma)\eqno(1.11)$$ Indeed $\ell(\sigma)$ walls of the form $\hh_{\varphi,n}$,
where $n\neq 0$, clearly separate
$Reg(w(A_{\sigma}))$ from $Reg(w(A_1))$ but $\ell(w)$ walls of the form $\hh_{\varphi,0}$ separate  
$Reg(w(A_1))$ from $Reg(A_1)$.  \vs Recall that $\psi$ is the highest root. For any 
integer $k\in \Bbb Z_+$ let
$\hh^{(k)} = \{x\in \hh^+\mid (\psi,x)\leq k\}$. Clearly $\hh^{(k)} =
k\,A_1$ so that, if $k\in \Bbb N$, the interior of $\hh^{(k)}$ is given by $$Int(\hh^{(k)}) =
\{x\in\hh\mid (\alpha_i,x)>0,\,\forall i\in
I\,\,\hbox{and}\,(\psi,x)<k\}\eqno (1.12)$$ Obviously 
$\hh^{(k)}$ is the closure of its interior. Furthermore since every point
in
$\hh$ lies in at most a finite number of alcoves (e.g. from volume
considerations) it follows from Remark 1.2 where
$\varphi = \psi$ and $\varphi=\alpha_i,\,i\in I$, that, if $k\in \Bbb N$,
$\hh^{(k)}$ is a union of alcoves. The following very simple and neat
observation (and proof) was made in [Cellini and Papi]. \vskip 1pc {\bf
Proposition 1.4.} {\it There are exactly $k^{\ell}$ alcoves in
$\hh^{(k)}$.} \vs {\bf Proof} (Cellini and Papi). Since $\hh^{(k)} = k\,A_1$ the
 volume of $\hh^{(k)}$ is $k^{\ell}$ times the volume of $A_1$.
The result then follows since every alcove necessarily has the same
volume. QED\vskip 1pc {\bf Remark 1.5.} Assume $\sigma\in W_f^+$.
If $\ell(\sigma)<k$ note that $$A_{\sigma}\s \hh^{(k)}\eqno (1.13)$$ Indeed otherwise
the $\psi$-walls $\hh_{\psi,j},\,j=1,\ldots,k$, would
separate $Reg(A_{\sigma})$ from $Reg(A_1)$ contradicting the fact that
$\ell(\sigma)<k$. QED \vskip 1pc It follows from Proposition 1.4 and Remark 1.5 that there
exists a formal power series $P(t)$ with coefficients in $\Bbb
Z_+$ such that $$P(t) = \sum_{\sigma\in W_f^+}
t^{\ell(\sigma)}\eqno (1.14)$$ The alcoves in $\hh^+$ have a
well-known connection with the loop group $\Omega(K)$. See [B]
and p.444 in [Ku-2]. Conforming to much of current terminology we
take the  exponents $\{m_i\},\,i\in I$, of $K$ to be the positive
integers (in non-decreasing order) such that the product of $(1+
t^{2m_i +1})$ over $i\in I$ is the Poincar\'e polynomial of $K$.
This makes 
$m_i$ here have value 1 less than the value assigned to $m_i$ in
[B]. The following is a classic result of Bott on the Poincar\'e
series of $\Omega(K)$. See Theorem B and (13.2) in [B]. \vskip 1pc {\bf Theorem 1.6} (Bott). 
{\it The Poincar\'e series of
the loop group $\Omega(K)$ is $P(t^2)$. Furthermore $$P(t) =
\prod_{i\in I} 1/(1-t^{m_i})\eqno (1.15)$$} \vskip .5pc 1.5. Let
$D\s \hh^+$ be the set of all dominant integral linear forms on $\hh$ and
for each $\lambda\in D$ let $\pi_{\lambda}:K\to Aut\, V_{\lambda}$ be
an irreducible representation with highest weight $\lambda$. As
usual $\pi_{\lambda}$ will also denote the corresponding
representation of $\g$ and the universal enveloping
algebra
$U(\g)$ on $V_{\lambda}$. Let $Cas\in Cent(U(\g))$
be the quadratic Casimir element corresponding to the Killing
form. For any $\lambda\in D$ let $Cas(\lambda)$ be the value of
the infinitesmal character of $\pi_{\lambda}$ on $Cas$. We
recall that $Cas(\lambda) = (\lambda+ \rho,\lambda+ \rho) -
(\rho,\rho)$ or $$Cas(\lambda) = (\lambda,\lambda) +
(2\,\rho,\lambda) \eqno(1.16)$$ where as usual $\rho =
1/2\,\sum_{\varphi\in \Delta_+} \varphi$. 

Recall that $\psi$ is the
highest root. As in \S 2.2
in [Ko-3] let
$$h_P = 1/(\psi,\psi)\eqno (1.17)$$ Of course 
$\psi$ is the highest weight of the adjoint representation so that
$Cas(\psi) =1$. Thus as already noted in  (2.2.3) in [Ko-3] one has
(see (1.16)) $1 =(2\rho,\psi) + (\psi,\psi)$. This immediately
implies that $$h_P = (2\rho,\psi)/(\psi,\psi) +1
\eqno  (1.18 )$$ (see (2.2.4) in [Ko-3]). Since $\psi$ is a
long root one has $h_P = 1/(\varphi,\varphi)$ for any long root and
$h_P$ is a positive integer. Let $h$ be the Coxeter number of
$\g$. 

In the later publication, [Ka], the number $h_P$ was referred as to as the dual Coxeter
number and was denoted by $h^{\vee}$. It plays a major role in Kac-Moody theory. Conforming to this
now well accepted terminology one has
\vskip 1pc {\bf Proposition 1.8.} {\it  One has 
$1/(\varphi,\varphi)$ is the dual Coxeter number $h^{\vee}$ for any long
root $\varphi\in \Delta$. Furthermore $h^{\vee} = h$ if $\g$ is simply-laced.}\vs {\bf Proof.} For
the first statement see the argument at the end of exercise 6.2 in \S 6.8 of [Ka]. The last
statement is Proposition 2.2 in [Ko-3]. QED \vskip 1pc {\bf Remark 1.9.} One has
$1/(\varphi,\varphi)$ is a positive integral multiple of the dual
Coxeter number for any $\varphi\in
\Delta$. See Proposition 2.3.1 and its proof in [Ko-3]. In particular 
$1/(\varphi,\varphi)$ is a positive integer for any $\varphi\in
\Delta$. \vs 1.6. By definition (see \S1.1) for any $\sigma\in W_f$ there uniquely exists $w^{\sigma}
\in
W$ and $z^{\sigma}\in \Gamma$ such that for any $x\in \hh$ one has $$\sigma(x) = w^{\sigma}(x) +
z^{\sigma}\eqno (1.19)$$ Let $\Delta_- = - \Delta_+$ and for any $w\in W$ let $\Phi_w =
w(\Delta_-)\cap \Delta_+$ so that, as one knows, $$\ell(w) = card\,\Phi_w\eqno (1.20)$$
\vskip 1pc {\bf Proposition 1.10.} {\it Let $\sigma\in W_f^+$. Then $$\ell(\sigma) +
\ell(w^{\sigma}) = (2\rho,z^{\sigma})\eqno (1.21)$$} \vs {\bf Proof.} Let $x\in Reg(A_1)$ so that
$y\in Reg(A_{\sigma})$ where $y = \sigma(x)$. Let $\varphi\in \Delta_+$. By (1.7) and (1.8) one
has that
$$n_{\varphi}(\sigma) +1 > \varphi(y) > n_{\varphi}(\sigma)\eqno (1.22)$$ For notational convenience
put $w=w^{\sigma}$. But since
$\varphi(w(x)) = w^{-1}(\varphi)(x)$ one has, by (1.19), $$\varphi(y) = w^{-1}(\varphi)(x)
+ \varphi(z^{\sigma}) \eqno (1.23)$$ But $0<|w^{-1}(\varphi)(x)|<1$. Since $\varphi(z^{\sigma})$ is an
integer, one has, by (1.22), $\varphi(z^{\sigma})= n_{\varphi}(\sigma)$ or $n_{\varphi}(\sigma) +1$,
according as
$w^{-1}(\varphi)\in \Delta_+$ or $w^{-1}(\varphi)\in \Delta_-$, that is, according as $\varphi\notin
\Phi_w$ or $\varphi\in \Phi_w$. But then summing $\varphi(z^{\sigma})$ over all $\varphi\in
\Delta_+$ yields (1.21), by (1.6) and (1.20). QED\vskip 1.5pc \centerline{\bf 2. Powers of the Euler
product and the set of weights $D_{alcove}$}\vskip 1.5pc 2.1.  Let
$i\in \I$.  One knows $(\rho,\alpha_i^{\vee}) = 1$ so that
$$(2\rho,\alpha_i) = (\alpha_i,\alpha_i)\eqno (2.1)$$ On the other
hand by (1.22) and (1.23) one has $$(2\rho,\psi) = 1-
(\psi,\psi)\eqno (2.2)$$ It follows from (2.1) and (2.2) that
$$2\rho \in Reg(A_1)\eqno (2.3)$$ As in [Ko-3] (see \S 3.1) let
$a_P\in K$ be defined by putting $a_P= exp\,2\pi i\,\,2\rho$. Let
$a\in K$. In [Ko-3], \S 3.1, we said that $a$ will be called an
element of type $\rho$ if it is conjugate to $a_P$. (Because of
the factor 2 this choice of terminology is perhaps inappropriate but it will be
retained nontheless.)
 For any $\lambda\in \Lambda$ let
$\chi_{\lambda}$ be the $K$-character of the irreducible
$\pi_{\lambda}$. In [Ko-3],
\S 3.1, we proved the following theorem about the $\hbox{dim} \,K$ power of the
Euler product $\prod_{n=1}^{\infty} (1-x^n)$. 
\vskip 1pc {\bf Theorem 2.1.} {\it For any
$\lambda\in D$ one has $\chi_{\lambda}(a_P)\in \{-1,0,1\}$ and as
formal power series $$(\prod_{n=1}^{\infty} (1-x^n))^{dim\,K} =
\sum_{\lambda\in
D}\chi_{\lambda}(a_P)\,\,dim\,V_{\lambda}\,\,x^{Cas(\lambda)}\eqno
(2.4)$$} \vskip 1pc See Theorem 3.1 in [Ko-3]. \vskip
1pc{\bf Remark 2.2.} For a relation between (2.4) and the Laplace
-Beltrami operator on $K$ (and implicitly the heat kernel on $K$)
see \S 5 in [Ko-3] and [F]. See also [Ze] for a recent physical application of the
heat kernel aspect of (2.4).
\vs 2.2. For any
$\sigma\in W_f^+$ let
$\lambda^{\sigma}\in \hh$ be defined by putting $$\lambda^{\sigma} = (w^{\sigma}(\rho)-\rho) +
z^{\sigma}/2\eqno (2.5)$$\vskip .5pc {\bf Proposition 2.3.} {\it Let $\sigma\in W_f^+$. Then
$\lambda^{\sigma}\in D$. Furthermore $\lambda^{\sigma}$ is in the root lattice.} \vs {\bf Proof.} 
Since $w^{\sigma}(\rho)$ and $\rho$ are two weights of
the representation $\pi_{\rho}$ it follows that $w^{\sigma}(\rho)
-\rho$ is in the root lattice. But, by (1.3), $z^{\sigma}/2$ is in the lattice spanned
by $\alpha_i/(\alpha_i,\alpha_i),\,i\in I$. Hence $z^{\sigma}/2$ is in the root lattice by Remark 1.9.
Thus $\lambda^{\sigma}$ is in the root lattice. In particular $\lambda^{\sigma}$ is in the weight
lattice. But then $\lambda^{\sigma}+ \rho $ is in the weight lattice. 
 But by the definition (2.5) one immediately has
$$2(\lambda^{\sigma}+ \rho) = \sigma(2\rho)\eqno (2.6)$$ Thus $$2(\lambda^{\sigma}+ \rho)\in
Reg(A_{\sigma})\eqno (2.7)$$ In particular $2(\lambda^{\sigma}+ \rho)\in Int(\hh^+)$. But then
$\lambda^{\sigma}+
\rho\in Int(\hh^+)$. But this implies that $\lambda^{\sigma}$ is dominant so that
$\lambda^{\sigma}\in D$. QED \vskip 1pc  Let $D_{alcove}=
\{\lambda^{\sigma}\mid
\sigma\in W_f^+\}$. By drawing attention to $D_{alcove}$ we have in
effect isolated what will be seen to be a distinguished subset of $D$
or, more significantly, a distinguished subset of the set of all
finite dimensional irreducible representations of $K$. It is
obvious from (2.7) that the map
$$W_f^+
\to D_{alcove}\,\qquad \sigma\mapsto \lambda^{\sigma}\eqno (2.8)$$ is bijective so that the subset
$D_{alcove}$ is parameterized by the set of all alcoves in a Weyl
chamber. \vskip 1pc 
Obviously in the formula (2.4) the only contributions to the sum on
the right hand side correspond to those $\lambda\in D$ such that
$\chi_{\lambda}(a_P) \in \{-1,1\}$. We now find that this condition
characterizes $D_{alcove}$. \vskip 1pc {\bf Theorem 2.4.} {\it Let
$\lambda\in D$ so that $\chi_{\lambda}(a_P) \in \{-1,0,1\}$. Then
$\chi_{\lambda}(a_P) \in \{-1,1\}$ if and only if $\lambda \in
D_{alcove}$. In particular the equality (2.4) simplifies to 
$$(\prod_{n=1}^{\infty} (1-x^n))^{dim\,K} =
\sum_{\sigma\in
W_f^+}\chi_{\lambda^{\sigma}}(a_P)\,\,dim\,V_{\lambda^{\sigma}}\,\,x^{Cas(\lambda^{\sigma})}\eqno (2.9
)$$ Furthermore if $\sigma\in W_f^+$ then 
$$\eqalign{\chi_{\lambda^{\sigma}}(a_P) &= (-1)^{\ell(\sigma)}\cr &= (-1)^{\ell(w^{\sigma})}\cr}\eqno
(2.10)$$}\vskip 1pc {\bf Proof.} We will use results and
notations of [Ko-3]. By Lemma 3.6 in {\it loc. cit.} one has
$$\chi_{\lambda}(a_P) = \epsilon_P(\lambda)\eqno (2.11)$$ where
$\epsilon_P(\lambda)$ is defined following Lemma 3.5.2, p. 199. 
The definition of $\epsilon_P(\lambda)$ rests upon an earlier
definition of the lattice
$M_P$. This lattice is defined in the line preceding
Proposition 2.3.1. In the notation of the present paper $M_P =
\Gamma/2$. Lemma 3.5.2 can then be restated as follows: either (1)
$v(\lambda+\rho)-\rho\notin \Gamma/2$ for all $v\in W$ or (2) there
exists a unique $v\in W$ such that $v(\lambda+\rho)-\rho\in
\Gamma/2$. 
Lemma 3.5.2 can then be restated as follows: either (1)
$v(\lambda+\rho)-\rho\notin \Gamma/2$ for all $v\in W$ or (2) there
exists a unique $v\in W$ such that $v(\lambda+\rho)-\rho\in
\Gamma/2$. By definition $\epsilon_P(\lambda) = 0$ in case (1) and 
$\epsilon_P(\lambda) = (-1)^{\ell(v)}$ in case (2). But then,
by (2.11),
$\chi_{\lambda}(a_P)=0$ in case (1) and $$\chi_{\lambda}(a_P)=
(-1)^{\ell(v)}\eqno (2.12)$$ in case (2). We will prove that
$\lambda\in D_{alcove}$ if and only if $\lambda$ satisfies the
condition of case (2). Assume case (2). Let
$v\in W$ be such that $v(\lambda+\rho)-\rho\in \Gamma/2$. Then there
exists
$\gamma\in \Gamma$ such that $$v(2(\lambda+\rho))-2\rho = \gamma$$
Applying $v^{-1}$ to this equality yields $$2(\lambda+\rho) =
w(2\rho) + z\eqno (2.13)$$ where we have put $w = v^{-1}$ and $z=
w(\gamma)$. But clearly $z\in \Gamma$ since $\Gamma$ is stable
under the action of $W$. Let $\sigma\in W_f$ be defined by putting
$\sigma = t_{z}\,w$. Then (2.13) asserts that $$2(\lambda+\rho) =
\sigma(2\,\rho)$$ But $2(\lambda+\rho)\in Int(\hh^+)$ since $\lambda\in D$. Thus $\sigma\in W_f^+$ and
hence 
$\lambda= \lambda^{\sigma}\in D_{alcove}$ by (2.6). Also $z= z^{\sigma}$ and
$w= w^{\sigma}$. But now $\ell (v) = \ell(w^{\sigma})$ since
$w^{\sigma} = v^{-1}$. But now clearly
$(2\rho,2\,\varphi/(\varphi,\varphi)\in 2\,\Bbb Z$ for any
$\varphi\in \Delta$. Thus $(2\,\rho,z^{\sigma})\in 2\,\Bbb Z$.
But then, recalling (1.21), this implies that the parity
of $\ell(\sigma)$ is the same as the parity of $\ell(w^{\sigma})$. Consequently 
$$\chi_{\lambda^{\sigma}}(a_P) =
(-1)^{\ell(\sigma)}\eqno (2.14)$$ by (2.12).

Now conversely assume $\lambda\in D_{alcove}$ so that
$\lambda=\lambda^{\sigma}$ for a unique $\sigma\in W_f^+$. Thus
$2\,(\lambda^{\sigma} + \rho) = w^{\sigma}(2\,\rho) + z^{\sigma}$ by (2.6).
Applying $v$ where $v^{-1}= w^{\sigma}$ and dividing by 2 
yields the relation $v(\lambda + \rho) = \rho + \gamma$ where
$\gamma = v(z^{\sigma})/2$. But then $\gamma\in \Gamma/2$ and hence
$\lambda$ satisfies the condition of case (2). The result then
follows from (2.12) and (2.14). QED\vs 2.3. 
In this section we will consider the
case where $K = SU(m)$ (i.e
$\g_{\Bbb C}$ is of type $A_{m-1}$. In this case combinatorists have come upon the set $D_{alcove}$ from
an entirely different perspective. Let $Q$ be the set of all partitions $q = (q_1,\ldots,q_m)$ of
length at most $m-1$ so that $q_m= 0$. One has a bijection $$f:D\to Q \eqno (2.15) $$ where if
$q = f(\lambda)$ then $(\lambda,\alpha_i^{\vee})= q_i-q_{i+1},\,i\in I$, where
$e_{\alpha_i}$ is the matrix unit $e_{i,i+1}$. 

Associated to $q\in Q$ is another partition $\widetilde {q}$ called $m$-core of $q$. See Exercise
1.1.8(c) in [Ma-2] or \S 3.4 in [A-F] or p. 467-469 in [St]. The partition
$\widetilde {q}$ is derived from $q$ by a step-by-step process of appropriately removing, from the
Young diagram of $q$, what are called
$m$-border strips in [Ma-2], [St] or rim hooks of length $m$ in [A-F]. The process is terminated when
no more removals are possible. What remains is the Young diagram of $\widetilde {q}$. The proof that
$\widetilde {q}$ is uniquely determined is particularly nice in [A-F]. One says $q$ has a null $m$ core
if
$\widetilde {q}$ is the empty partition. Let $Q_o = \{q\in Q\mid q\,\hbox{has null $m$
core}\}$. It is obvious that the size of any $q\in Q$ is a multiple of $m$. See Lemma 3.4 in
[A-F] for a neat characterization of the elements in $Q_o$. \vskip 1pc {\bf Theorem 2.5. } {\it
One has $f(D_{alcove})\s Q_o$ and $$f:D_{alcove}\to Q_o\eqno (2.16) $$ is a bijection.}
\vskip 1pc {\bf Proof.} Let $\lambda\in D$. Then by 1.3.17(a), p. 50 in [M] one has
$\chi_{\lambda}(a_P) \in \{0,1,-1\}$ and $\chi_{\lambda}(a_P) \in \{1,-1\}$ if and only if $\lambda\in
Q_o$. A similar statement is made in Theorem 5.7 of [A-D]. (Of course these are statements for
the
$SU(m)$ case.) The result then follows from Theorem 2.4 in the present paper. QED \vskip 1pc {\bf Remark
2.6.} 
Using Bott's formula (1.13),
 one can show that the number $q\in Q_o$ having size
$mk$ is $m + k -2\choose m-2$. However both R. Stanley and R. Adin have pointed out that this
statement can be deduced from known facts about $Q_o$.  
\vskip 1.5pc

\centerline{\bf 3. The structure of the homology $H_*(\u^-)$.}\vskip 1.5pc

3.1. The main
results of this paper are given in
\S4. We have been convinced by Pavel Etingof that the results are best illuminated using results of
Garland [G], Garland-Lepowsky [G-L] and [Ku-1,2] on the homology and cohomology of the ``nilradical" of
the standard maximal parabolic subalgebra of the affine Kac-Moody Lie algebra associated to $\g$. For the
definition of the affine Kac-Moody Lie algebra see \S 6 in [Ka] or \S 13.1 in [Ku-2]. In the main we
will adhere to the development of affine  Kac-Moody Lie algebras in these references except for some
changes in notation. In particular
 we retain our previous meaning of
$\g$ (and not write 
$\buildrel\circ \over\g$ for the finite dimensional complex simple Lie algebra). Also we write $\www
{\g}$ for the affine Kac-Moody Lie algebra $\hat{{\cal L}}(\buildrel\circ \over\g)$ given in (1), p. 482
in [Ku-2]. One has $\g\s \www\g$ and if $c,d\in\www \g$ are as in \S 13.1.1 in [Ku-2] then, writing
$\delta$ for $c$, the complexification $\www \hh_{\Bbb C}$ of the real $\ell + 2$-dimensional abelian Lie
subalgebra
$$\www
\hh =
\hh +
\Bbb R \delta +
\Bbb R d\eqno (3.1)$$ ``serves" as  a Cartan subalgebra of $\www \g$. Let $\d = \Bbb R \delta +
\Bbb R d$. We extend the positive definite bilinear form $B|\hh$ to a nonsingular bilinear form $(x,y)$ on
$\www \hh$ so that $\hh$ is orthogonal to $\d$, $(d,d) = (\delta,\delta) = 0$ and $(d,\delta) = 1$. This
is further extended to $\www \hh_{\Bbb C}$ by complex linearity. Using the latter extension we 
identify $\www \hh_{\Bbb C}$ with its dual space. As a linear space one has the direct sum
decomposition $$\www\g = \g_{\d} +  t\,\g[t] +  t^{-1}\,\g[t^{-1}] \eqno (3.2)$$ where
$\d_{\Bbb C}$ is the complexification of
$\d$ and we have put $\g_{\d} = \g +\d_{\Bbb C}$. For a parameter $u$ the space $\g[u]$ is the direct sum
$$\g[u] = \sum_{k=0}^{\infty} u^k\, \g\eqno (3.3)$$ The commutation relations in $\www \g$ are given
in (2), \S 13.1.1 in [Ku-2]. In particular $\d_{\Bbb C} = Cent\,\g_{\d}$ and the adjoint action of
$\g$ on $\www\g$  stabilizes $\Bbb C[t,t^{-1}]\g$ and 
 is linear with respect to the $\Bbb C[t,t^{-1}]$-module structure on $\Bbb C[t,t^{-1}]\g$. It is
otherwise obvious. One has $ad\,\delta = 0$ and $ad\,d $ is the $t$-degree operator $t\,{d \over dt}$. 

For any $\varphi\in \Delta$ let $0\neq e_{\varphi}\in \g$ be a
corresponding weight vector. The set
$\www{\Delta}'$ of affine roots is just the set of nonzero weights for the adjoint action of
$\www\hh_{\Bbb C}$ on $\www \g$. For any affine root $\beta$ let $\www \g_{\beta}\s
\www\g$ be the corresponding root space. One has the disjoint union $$\www{\Delta}' =
\www{\Delta}^{Im}\cup
\www{\Delta}$$ where $\www{\Delta}^{Im}= (\Bbb Z-\{0\})\delta$. The elements in  $\www{\Delta}^{Im}$ are
called imaginary affine roots. If
$\beta$ is an imaginary affine root so that
$\beta = m\,\delta$ for a nonzero integer $m$, then $\www \g_{\beta}= t^m\,\hh_{\Bbb C}$. One has $$
\www{\Delta}= \Delta + \Bbb Z\,\delta\eqno (3.4)$$ One notes that $$\www{\Delta}'\s
\www\hh\eqno (3.5)$$ The elements of
$\www{\Delta}$ are called real affine roots. If $\beta$ is a real affine root, so that $\beta = \varphi +
k\,\delta$, where
$\varphi\in
\Delta$ and $k\in \Bbb Z$, put $$e_{\beta} =
t^k\,e_{\varphi}\eqno (3.6)$$ (this is clearly unambiguous even if $\beta = \varphi$) and one has
$$\www\g_{\beta} =\Bbb Ce_{\beta}\eqno (3.7)$$  One introduces the set $\www{\Delta}'_+$ of positive
affine roots by putting  $$\www{\Delta}'_+= \www{\Delta}^{Im}_+\cup
\www{\Delta}_+\eqno (3.8)$$ where $\www{\Delta}^{Im}_+ = \Bbb N\delta$ and $$\www{\Delta}_+ =
(\Delta_+ +
\Bbb Z_+\,\delta) \cup (\Delta_- + \Bbb N\,\delta)\eqno (3.9)$$ The sets
obtained by multiplying each of the 3 sets in (3.8) by $-1$ is denoted by replacing the subscript $+$ by
the subscript $-$. If $\beta$ is a real affine root then clearly $(\beta,\beta)>0$ (see (3.4) and (3.5))
so that $\www\hh_{\beta}$, defined as the orthocomplement of $\beta$ in $\www\hh$, is a subspace of
codimension 1 in $\www\hh_{\beta}$. Let $s_{\beta}$ be the (linear) orthogonal reflection of $\www\hh$
defined by the subspace $\www\hh_{\beta}$. Let $I_*= I\cup \{0\}$ and let $\alpha_0\in \www{\Delta}_+$
be defined by putting $\alpha_0=\delta-\psi$. The ``Weyl group" of $\www\g$ is the Coxeter group $\www
W$ with simple generators $\{s_{\alpha_i}\},i\in I_*$, operating linearly on $\www\hh$. See \S 1.3.1
in [Ku-2]. As such one has a length function $\tau\mapsto \ell(\tau)$ on $\www W$. Also
 $s_{\beta}\in \www W$ for any real affine root $\beta$. Since the simple generators of $\www W$
include the simple generators of $W$ one has a natural embedding of
$W$ in $\www W$. The definition of $\Phi_w$ for $w\in W$ (see \S1.6) extends to any $\tau \in\www W$ by
putting
$\Phi_{\tau} = \www{\Delta}'_+ \cap \tau(\www{\Delta}'_-)$. Since $\delta$ is fixed under the action of
$\www W$ one notes that $$\Phi_{\tau} = \www{\Delta}_+ \cap \tau(\www{\Delta}_-)\eqno (3.10)$$
Furthermore not only is $\Phi_{\tau}$ a finite set but in fact $$\ell(\tau) = card\,\Phi_{\tau}\eqno
(3.11)$$ See e.g. Lemma 1.3.14 in [Ku-2]. 

Let $\www\hh_1$ be the hyperplane in
$\www\hh$ defined by putting $\www\hh_1 =\{x\in\www\hh\mid (\delta,x) = 1\}$ so that $\www\hh_1 = d + \Bbb
R\,\delta + \hh$. It is clear that $\www\hh_1$ is stable under the action of $\www W$. Let $\zeta:\www\hh
\to
\hh$ be the projection with kernel $\d$. Then
$$\www W\to W_f,\qquad
\tau
\mapsto \overline \tau\eqno (3.12) $$ is a group isomorphism where for any $x\in \hh$, $$\overline
{\tau}(x) = \zeta(\tau (d+x))\eqno (3.13)$$ See \S 6.6 in [Ka]. If $\varphi\in \Delta_+,\,n\in \Bbb N$
and and $k\in \Bbb Z_+$ one readily notes that $$\eqalign{\overline {s_{n\delta - \varphi}} &=
s_{\varphi,n}\cr  \overline {s_{k\delta + \varphi}}&= s_{\varphi,-k}\cr}\eqno (3.14)$$ See e.g. p. 132
in [C-P]. Extend the map (3.12) to be an involutory bijection on the set $\www W\cup W$ by putting
$$\overline {\sigma}= \tau\eqno (3.15)$$ where $\sigma\in W_f$ and $\tau\in \www W$ is such that
$\overline {\tau} =
\sigma$. 

Let $\sigma\in W_f^+$. It is clear that any wall which separates $Reg(A_1)$ from $Reg(A_{\sigma})$ is
necessarily of the form $\hh_{\varphi,n}$ where $n>0$ and $\varphi\in \Delta_+$. On the other hand it is
immediate, say, from Lemma 1.3.14 in [Ku-2] and (3.14), (See also (1.1), p. 132 in [C-P]) that
$\hh_{\varphi,n}$ is such a wall if and only if $n\delta -\varphi\in \Phi_{\overline {\sigma}}$. That is,
one has \vs {\bf Proposition 3.1.} {\it Let
$\sigma \in W_f^+$. Then $$\Phi_{\overline {\sigma}} = \{n\delta-\varphi\mid
\hh_{\varphi,n}\,\,\hbox{separates $Reg(A_1)$ from $Reg(A_{\sigma})$}\}$$}\vskip .5pc 3.2.
Let $$\www\rho = d/2 + \rho \eqno (3.16)$$ We note that for $i\in I_*$ one has
$$(\www\rho,\alpha_i) = (\alpha_i,\alpha_i)/2\eqno (3.17)$$ Indeed the equality (3.17) is well known
for
$i\in I$. Since $\alpha_0 = \delta-\psi$, for $i = 0$, one has $$\eqalign{(\www\rho,\alpha_0)&= 1/2
- (\rho,\psi)\cr &= (\psi,\psi)/2\cr &= (\alpha_0,\alpha_0)\cr}$$ by (1.17) and (1.18). This proves
(3.17).

 For any subset $\Phi\s \www {\Delta}_+$
let $\langle\Phi\rangle = \sum_{\beta\in \Phi} \beta $ . For the proof of the following
((3.18)) extension of (5.10.1) in [Ko-1] to the Kac-Moody case see Proposition 2.5 in [G-L] or (3) in
Corollary 1.3.22 in [Ku-2]. For any $\tau\in \www W$ one has $$\www\rho- \tau(\www\rho) =
\langle\Phi_{\tau}\rangle \eqno (3.18)$$ For any $\sigma\in W_f$ and $\varphi\in \Delta_+$ recall
the definition, in (1.5), of the integer $n_{\varphi}(\sigma)$.  \vs 
 {\bf Lemma 3.2.} {\it Let $\sigma\in W_f^+$. Then $$\www\rho-\overline {\sigma}(\www\rho) =
(\sum_{\varphi\in \Delta_+}n_{\varphi}(\sigma)(n_{\varphi}(\sigma)+1)/2)\,\delta - 
(\sum_{\varphi\in \Delta_+}n_{\varphi}(\sigma)\,\varphi)\eqno (3.19) $$}\vs {\bf Proof.} Let $\varphi
\in \Delta_+$. It is immediate from (1.7) and (1.8) that the wall $\hh_{\varphi,j}$ separates $Reg(A_1)$
from $Reg(A_{\sigma})$ if and only if $j$ is a positve integer such that $1\leq j\leq
n_{\varphi}(\sigma)$. But then (3.19) follows immediately from Proposition 3.1 and (3.18) where $\tau =
\overline {\sigma}$. QED \vskip 1pc One now has the following explicit expression for the elements in
$D_{alcove}$ (see \S 2.2). \vs {\bf Theorem 3.3. } {\it Let
$\sigma\in W_f^+$ so that
$\lambda^{\sigma}\in D_{alcove}$ (see \S 2.2). Then $$\lambda^{\sigma} = \sum_{\varphi\in
\Delta_+}n_{\varphi}(\sigma)\,\varphi\eqno (3.20)$$}\vs {\bf Proof.} By definition $\www\rho = d/2 + \rho$
(see (3.16)). Thus $2\,\www\rho = d + 2\rho$. On the other hand by the definition of $\overline{\sigma}$
(see (3.12), (3.13) and (3.15)) one has $$\eqalign{\sigma(2\,\rho) &= \zeta(\overline {\sigma}(d +
2\rho)\cr &=p(\overline {\sigma}(2\,\www\rho)} \eqno (3.21)$$ But $2\,\rho = \zeta(2\,\www\rho)$. But then
$$\sigma(2\,\rho) - 2\rho = \zeta(\overline {\sigma}(2\,\www\rho) - \zeta(2\,\www\rho)$$ But
$\overline {\sigma}$ is linear. Thus, by
(3.19),
$$\eqalign{\sigma(2\,\rho) - 2\rho &= 2\,\zeta(\overline {\sigma}(\www\rho) - \www\rho)\cr
&= 2\sum_{\varphi\in \Delta_+}n_{\varphi}(\sigma)\,\varphi)\cr}\eqno (3.22)$$ But 
$\sigma(2\,\rho) - 2\rho
= 2\,\lambda^{\sigma}$ by (2.6). But then (3.19) follows from (3.22). QED \vs {\bf Remark 3.4.}. As
mentioned in the proof above $\overline {\sigma}$ is linear and we have used this fact. However $\sigma$
is not linear and (2.6) does not imply that $\sigma(\rho) - \rho = \lambda^{\sigma}$. In fact in general 
$\sigma(\rho) - \rho \neq \lambda^{\sigma}$. Instead one has $\sigma(2\rho)/2 - \rho =
\lambda^{\sigma}$ by (2.6).\vs

\def\wd{\wedge\u^-}

\rm 
3.3. Recall (3.2). For notational simplicity put $\u = t\,\g[t]$ and $\u^{-}=
t^{-1}\,\g[t^{-1}]$. To state the results of [G], [G-L] and [Ku-1,2] on the homology $H_*(\u^-)$
it is clarifying to write down a $\www \hh$-weight basis of the
exterior algebra $\wedge \u^-$. 

Let $J = \{1,2,\ldots,dim\,\k\}$ and let $x_j,\,j\in J$, be a
$\hh$-weight basis of $\g$ (under the adjoint representation). 
For any $j\in J$ let $\mu_j\in \hh$ be the weight corresponding to
$x_j$ so that $\mu_j\in \Delta\cup \{0\}$.  For any $n\in \Bbb N$
 let $I_n = \{1,\ldots,n\}$ and $P_n$ be the set of all partitions
$p = (p_1,\ldots,p_n)$ of length (exactly) $n$. Let
$R_n$ be the set of all maps
$$r:I_n\to \Bbb N\times J\eqno (3.23)$$ where if $r(i) =
(r_i,r_{[i]})$ then (1) $p(r)  = (r_1,\ldots,r_n)\in P_n$ and (2) if
$i<j$ and $r_i = r_j$ then $$r_{[i]} > r_{[j]}\eqno (3.24)$$ Given
$r\in R_n$ let $$z_r = t^{-r_1}\,x_{r_{[1]}}\wedge \cdots \wedge
t^{-r_n}\,x_{r_{[n]}}\eqno (3.25)$$ so that $$z_r\in
\wedge^n\u^-$$ Let
$\mu(r) =
\sum_{i\in I_n}
\mu_{r_{[i]}}$. The size,
$|p_r|$, of the partition $p_r$ equals $\sum_{i\in I_n}r_i$. Since
$d\in
\www\hh$ operates as $t{d\over dt}$ one immediately has \vskip 1pc {\bf
Proposition 3.5.} {\it Let $r\in R_n$. Then $z_r\in \wedge^n\u^-$ is
a $\www \hh$-weight vector of weight $$-|p(r)|\delta + \mu(r)$$}\vs
The condition (3.24) in the definition of $r$ guarantees that
the elements $z_r,\,r\in R_n$, are linearly independent. In fact
one immediately notes
\vs {\bf Proposition 3.6.} {\it The set $\{z_r\},\,r\in R_n$, is a
$\www\hh$-weight basis $\wedge^n\u^-$.} \vs Let $R = \cup_{n\in \Bbb Z_+} R_n$ so that 
$\{z_r\},\,r\in R$, is a
$\www\hh$-weight basis $\wedge\u^-$. If $Y$ is any $\www\hh$-module and $k\in
\Bbb Z$ we will denote the eigensubspace of $Y$, with eigenvalue (t-weight) $-k$, for the action of
$d$, by $(Y)_k$. Clearly $\{z_r\},\,r\in R,\,|p(r)|= k,$ is a basis of $(\wedge\u^-)_k$. It is
immediate from (3.25) that $(\wedge\u^-)_k$ is finite dimensional and one has the direct sum 
$$\wedge\u^- =
\sum_{k\in \Bbb Z_+}(\wedge\u^-)_k\eqno (3.26)$$ 

3.4. Now $\u^-$ (and of course also $\u$) is a Lie subalgebra
of $\www\g$ and $\wedge\u^-$ is a chain complex for the Lie algebra homology space
$H_*(\wedge\u^-)$. Since $\g_{\d}$ (see (3.2)) obviously normalizes $\u^-$ one notes that
$\wedge\u^-$ and $H_*(\wedge\u^-)$ are completely reducible $\g_{\d}$-modules with finite dimensional
irreducible components. In fact since
$d$ is central in
$\www\g$ it is immediate that
$(\wedge\u^-)_k$ is a finite dimensional, completely reducible, $\g$-module subcomplex (with
homogeneous components
$(\wedge^n\u^-)_k =
\wedge^n\u^-\cap (\wedge\u^-)_k)$ and $(H_*(\wedge\u^-))_k$ is just the homology of
$(\wedge\u^-)_k$. The decomposition (3.26) yields the direct sum, $\g$-module decomposition, with
completely reducible finite dimensional components,
$$H_*(\wedge\u^-) = \sum_{k\in
\Bbb Z_+}(H_*(\wedge\u^-))_k\eqno (3.27)$$ \vskip .5pc The determination of $H_*(\wedge\u^-)$ as a
$\g$-module is due to H. Garland. See Theorem 3.2 in [G]. This result was extended to the general
Kac-Moody case (and, in addition, with values in a suitable module) in [G-L]. See theorem 8.5 in
[G-L]. An elegant presentation of the Garland-Lepowsky theory is given in \S 3.2 of [Ku-2]. The
proof of Theorem 3.2 depends upon Theorem 2.5 in [G]. The latter (see Theorem 3.7 below), of
interest in itself, is a statement about the Laplacian operator $L$ (denoted by $\Delta$
in [G]) associated to the boundary operator of $\wedge\u^-$ and a positive definite
Hermitian structure on
$\wedge\u^-$. The
operator $L$ commutes with the action of
$\g_{\d}$ so that if $Har(\u^-) = Ker\,L$ and $Har_n(\u^-) = Ker\,L|\wedge^n\u^-$ then one has an
isomorphism $$Har(\u^-)\equiv H_*(\wedge\u^-)\eqno (3.28)$$ of graded $\g_{\d}$-modules. Thus it
suffices to explicitly determine $L$ and its kernel. This determination rests upon Theorem 2.5 of
[G]. As far as I am aware of, a proof of this theorem has not appeared in the literature. However the
result is established as a special case of a much stronger theorem (arbitrary symmetrizable
Kac-Moody case together with a suitable module) due to Kumar in [Ku-1]. Also see Theorem 3.4.2 in
[Ku-2] and the final remark on p. 107 in [Ku-2]. 

If $\xi\in \www\hh$ is an $\www\hh$-weight occurring in
$\wedge\u^-$ then $\xi = a\delta + \nu $ where $-a\in \Bbb Z_+$ and $\nu$ is in the ordinary root
lattice of
$\hh$. We may refer to
$a$ as the $\delta$ component of $\xi$ and $\nu$ as the $\hh$ component of $\xi$. The statement
that $\xi$ is dominant is just the statement that $\nu\in D$. In particular if $\xi$ is the
highest weight of a $\g_{\delta}$ irreducible component of $\wedge\u^-$ then certainly $\nu\in D$.
\vs {\bf Theorem 3.7} (Garland). {\it Let
$k\in
\Bbb Z_+$ and let $m_k$ be the maximal eigenvalue of $Cas$ on $(\wedge\u^-)_k$. Then $$m_k\leq k\eqno
(3.29)$$ Furthermore one has equality in (3.29) if and only if $(Har(\u^-))_k \not=  0$. Moreover in
such a case $(Har(\u^-))_k$ is the eigenspace in $(\wedge\u^-)_k$ for $Cas$ belonging to the
eigenvalue
$k$.} \vs {\bf Proof} (Kumar). Let $Z$ be an irreducible $\g_{\d}$-submodule of $(\wedge\u^-)_k$ and
let
$\xi$ be the highest weight. Then $\xi = -k\delta + \nu$ for some $\nu\in D$. Then $L|Z$ operates
as the scalar operator $$c = 1/2((\www\rho,\www\rho) - (\www\rho -k\delta + \nu,\www\rho -k\delta
+\nu))$$ by Theorem 3.4.2 in [Kum-2] since in the case at hand $\lambda = 0$. But $\www\rho = \rho
+ d/2$ by (3.16). Thus $c= 1/2((\rho,\rho) - (\rho + \nu,\rho + \nu) + k)$. That is $c= 1/2(k -
Cas(\nu))$. But $L|Z$ is positive semidefinite. This implies the inequality (3.29). But $c = 0$ if
and only if $Cas (\nu) = k$ and this must be the case if and only if $m_k = k$ and $Z\s
(Har(\u^-))_k$.  QED
\vs {\bf Remark 3.8.} Garland remarks that his Theorem 3.2 is an analogue of results in [Ko-1]. At
first glance Theorem 2.5 in [G], upon which his Theorem 3.2 depends, appears to have no analogue in
[Ko-1]. However Kumar's more general result, Theorem 3.4.2 in
[Ku-2], is in fact manifestly an infinite-dimensional analogue of Theorem 5.7 in [Ko-1]. 

The following statement is a corollary of Theorem 3.7. \vs {\bf
Theorem 3.9.} {\it Let $0\neq z$ be an $\www\hh$-weight vector in $\wedge\u^-$ and let $Z\s
\wedge\u^-$ be the
$\g_{\d}$-submodule generated by
$z$.  Let $\xi\in \www\hh$
be the weight of $z$ and let $\lambda$ be the $\hh$-component of $\xi$. Then
$Z\s Har(\wedge\u^-)$ and $z$ is a highest weight vector of $Z$ (so that $Z$ is an irreducible
$\g_{\d}$-submodule) if and only if the $\delta$-component of $\xi$ equals $(\rho,\rho)-
(\lambda +
\rho,\lambda + \rho)$. Moreover in such a case $\lambda\in D$ and $$\xi = -Cas(\lambda)\,\delta +
\lambda\eqno (3.30)$$} \vs {\bf Proof.} Let $-k$ be the $\delta$-component of $\xi$ so that $k\in
\Bbb N$ and
$Z\s (\wedge\u^-)_k$. Now if $Z\s Har(\wedge\u^-)$ then $Z\s (Har(\wedge\u^-))_k$ and hence $Cas|Z$
is the scalar operator for the scalar $k$, by Theorem 3.7. But if also $z$ is a highest weight
vector of
$Z$ then
$\lambda\in D$ and $$k =  (\lambda + \rho,\lambda +\rho) - (\rho,\rho)\eqno (3.31)$$ 

Conversely assume (3.31). If $z$ is not a highest weight vector of $Z$ there obviously exists
an irreducible
$\g_{\d}$-submodule
$Z'\s Z$, necessarily having $\xi$ as a weight, and such that if $\xi' = -k\,\delta + \lambda'$ is
the highest weight of $Z'$ then $\lambda \not=\lambda'$. But $\lambda'\in D$ and $Cas(\lambda') =
(\lambda' + \rho,\lambda' +\rho) - (\rho,\rho)$. But then $Cas(\lambda')> k$ by (3.31) and the
Freudenthal result (5.9.2) in [Ko-1]. This contradicts (3.29). Thus $z$ is a highest weight vector of
$Z$ and $\lambda\in D$. But then the right hand side of (3.31) equals $Cas(\lambda)$. Hence $Z\s
Har(\u^-)$ by Theorem 3.7. QED\vs 3.5. Theorem 3.2 in [G] gives
the decomposition of
$Har(\wedge\u^-)$ (and hence equivalently
$H_*(\wedge\u^-)$) as a
$\g_{\d}$-module. This result (see Theorem 3.10 below)
is also the application of Theorem 8.5 in
[G-L] to the case of the standard maximum parabolic subalgebra of the affine Kac-Moody Lie algebra
where the module is trivial. See also Theorem 3.2.7 in [Ku-2] for this case. The statement
is made stronger (implicit in [G-L]) by including Lemma 3.2.6 in [Ku-2]. One
notes that the strenthened statement is a Kac-Moody analogue of Lemma 5.12 and Theorem 5.14 in
[Ko-1]. By (3.18) one has
$$\overline{\sigma}(\www\rho) -
\www\rho = -\langle\Phi_{\overline
\sigma}\rangle\eqno (3.32)$$ for any $\sigma\in W_f$. 
\vs {\bf Theorem 3.10} (Garland). {\it Let $\sigma\in W_f^+$. Then $-\langle\Phi_{\overline
\sigma}\rangle
$ occurs as a $\hh_{\delta}$-weight of multiplicity one in the $\g_{\d}$- module $\wedge\u^-$. In
particular there exists a unique $r\in R$ (see \S3.3), henceforth denoted by $r^{\sigma}$, such that
$z_{r^{\sigma}}$ is a weight vector with weight $-\langle\Phi_{\overline \sigma}\rangle$. In the
notation of Lemma 3.5 $$-\langle\Phi_{\overline \sigma}\rangle = -|p(r^{\sigma})|\,\delta +
\mu(r^{\sigma})\eqno (3.33)$$ In addition $-\langle\Phi_{\overline \sigma}\rangle$ is dominant. That
is, $\mu(r^{\sigma})\in D$.

Let $Z_{\sigma}$ be the $\g_{\d}$-module generated by $z_{r^{\sigma}}$. Then $Z_{\sigma}$ is
$\g_{\d}$-irreducible and $z_{r^{\sigma}}$ is a highest weight vector of $Z_{\sigma}$. Moreover
$Z_{\sigma}\s Har(\wedge\u^-)$ and indeed one has the multiplicity free decomposition
$$Har(\wedge\u^-) =\sum _{\sigma\in W_f^+} Z_{\sigma}\eqno (3.34)$$ With respect to the two
compatible gradations, $Har_n(\wedge\u^-)$ and $(Har(\wedge\u^-))_k$ one has 
$$(Har(\wedge\u^-))_k = \sum_{\sigma\in W_f^+,\,\,\,|p(r^{\sigma})|= k}Z_{\sigma}\eqno (3.35)$$ and
$$Har_n(\wedge\u^-)= \sum _{\sigma\in W_f^+,\,\,\,\ell(\sigma)=n} Z_{\sigma}\eqno (3.36)$$}\vs {\bf
Remark 3.11.} Note that the existence of $r^{\sigma}$ with the cited properties is a consequence of 
Proposition 3.6 and the multiplicity one property of $-\langle\Phi_{\overline \sigma}\rangle$. Note
also that $(Har(\wedge\u^-))_k$ is finite dimensional since $(\wedge\u^-)_k$ is finite dimensional
(see \S 3.3). On the other hand $Har_n(\wedge\u^-)$ is also finite dimensional (even though
$\wedge^n\u$ is infinite dimensional) since the set $\{\sigma\in W_f\mid \ell(\sigma)= n\}$ is
obviously finite. \vs 3.6.  {\bf Remark 3.12.} Since the set
$E= \{t^{-i}\,x_j\mid (i,j)\in \Bbb N\times J\}$ is clearly a basis of $\u^-$ note that the elements
$z_{S_k},\,k\in I_n$, are linearly independent in $\wedge\u^-$ where $S_k,\,k\in I_n$, are finite 
mutually distinct subsets of $E$ and for any such subset $S$, $z_S$ is the decomposable element
obtained by exterior multiplication, in some order, of the elements in $S$.\vs 
Recalling the notation of (3.23) let $J^+ = \{j\in J\mid \mu_j \in
\Delta_+\}$. We may choose the ordering of the basis $\{x_j\}$ of $\g$ so that $J^+ =
\{1,\ldots,m\}$ where here $m = card\,\Delta_+$. For
$j\in J^+$ we now write $\varphi_j$ for
$\mu_j$ and choose $x_{j} = e_{\varphi_j}$. \vs For
$n\in
\Bbb Z_+$ let
$R_n^+$ be the set of all
$r\in R_n$ such that the image of (3.23) is contained in $\Bbb N\times J^+$. Now let $\sigma\in
W_f^+$ and let $n= \ell(\sigma)$. For $j\in J^+$ let $n_{\varphi_j}(\sigma)$ be defined as in
\S1.3 and let $z^{\sigma,j}\in \wedge^{n_{\varphi_j}(\sigma)}\u^-$ be defined by putting
$z^{\sigma,j}= 1$ if $n_{\varphi_j}(\sigma)= 0$ and otherwise
$$z^{\sigma,j} = t^{-1}\,x_j\wedge\cdots \wedge t^{-n_{\varphi_j}(\sigma)}\,x_j\eqno (3.37)$$ Next
let $z^{\sigma}\in \wedge^n\u^-$ (see (1.6)) be defined by putting $$z^{\sigma} =
z^{\sigma,1}\wedge \cdots\wedge z^{\sigma,m}\eqno (3.38)$$
\vskip.5pc We now use the results of \S 3.2 and relate Theorem 3.10 with $D_{alcove}$. See \S 2.2
and the cautionary Remark 3.4. \vs {\bf Theorem 3.13.} {\it Let $\sigma\in W_f^+$ and let
$\lambda^{\sigma}\in D_{alcove}$ be defined as in (2.5). Then in the notation of (3.33) one has
$\mu(r^{\sigma}) = \lambda^{\sigma}$ and in fact (3.33) can be written $$-\langle\Phi_{\overline
\sigma}\rangle = -Cas(\lambda^{\sigma})\,\delta +\lambda^{\sigma}\eqno (3.39)$$ In particular not
only is $Cas(\lambda^{\sigma})$ an integer but in fact $$Cas(\lambda^{\sigma}) = \sum_{\varphi\in
\Delta_+}n_{\varphi}(\sigma)(n_{\varphi}(\sigma)+1)/2\eqno (3.40)$$ Moreover (recalling
(3.38) and Theorem 3.10 ) one has $r^{\sigma}\in R_n^+$ and $$z_{r^{\sigma}} = \pm z^{\sigma}\eqno
(3.41)$$ Finally
$Z_{\sigma}\equiv V_{\lambda^{\sigma}}$ as a $\g$-module and, as a $\g$-module, $Har(\wedge\u^-)$ is
multiplicity free and one has the equivalence $$Har(\wedge\u^-) \equiv
\sum_{\sigma\in W_f^+} V_{\lambda^{\sigma}}\eqno (3.42)$$ Of course the same
statement is true when $H_*(\wedge\u^-)$ replaces $Har(\wedge\u^-)$.}\vs {\bf Proof.} The statement
that
$\mu(r^{\sigma}) = \lambda^{\sigma}$ is immediate from (3.19),(3.20),(3.32) and (3.33). Since
$z_{r^{\sigma}}$ and $Z_{\sigma}$ of Theorem 3.10 satisfies the condition of $z$ and $Z$ of Theorem
3.9 it follows from (3.30) that the $\delta$ component of $-\langle\Phi_{\overline \sigma}\rangle$
equals $-Cas(\lambda^{\sigma})$. This proves (3.39). But then (3.40) follows from (3.19). But now
the $\www\hh$ weight of the weight vector $z^{\sigma}$ is clearly $$(-\sum_{j=1}^m
n_{\varphi_j}(\sigma)(n_{\varphi_j}(\sigma) +1)/2)\,\delta + \sum_{j=1}^m
n_{\varphi_j}(\sigma)\varphi_j$$ Thus
$z^{\sigma}$ is an $\www\hh$ weight vector of weight $-\langle\Phi_{\overline \sigma}\rangle$ by
(3.19). Thus one has (3.41) by the multiplicity one condition (Theorem 3.10) of this weight in
$\wedge\u^-$. It follows in particular, (see Remark 3.12) that $r^{\sigma}\in R_n^+$. The remaining
statements are then obvious noting that the
$\hh$ component of
$-\langle\Phi_{\overline \sigma}\rangle$ determines the $\delta$-component (or using the
injectivity of the map (2.8)). QED\vs 3.7. Let $Q = (\Bbb Z_+)^m$ where $m =
card\,\Delta_+$. If $q\in Q$ let $q_i\in \Bbb Z_+,\,i\in I_m$, be defined so that $q=
(q_1,\ldots,q_m)$. Let $X\s \hh$ be the semigroup generated by $\Delta_+$ and let $\eta:Q\to X$
be defined by putting $$\eta(q) = \sum_{i\in I_m} q_i\,\varphi_i\eqno (3.43)$$ If $\eta(q) = \mu$
we will refer to $q$ as a positive root
partition of $\mu$. For any $\mu\in X$ let $$Q_{\mu} = \eta^{-1}(\mu)\eqno (3.44)$$ so
that $Q_{\mu}$ is the set of all positive root partitions of $\mu$. 

Now if $q\in Q$ let $$c(q) = \sum_{i\in I_m} q_i(q_i +1)/2\eqno (3.45)$$ and let $z^{(q)}\in
\wedge\u^-$ be defined by putting $$z^{(q)} = z^{(q),1}\wedge\cdots\wedge z^{(q),m}\eqno (3.46)$$
where
$z^{(q),i}=1$ if $q_i=0$ and otherwise $$z^{(q),i} = t^{-1}x_{i}\wedge \cdots \wedge
t^{-q_i}\,x_i$$ One notes that $z^{(q)}$ is an $\www\hh$-weight vector and the $\www\hh$
$$\hbox{weight of}\,\ z^{(q)}= -c(q)\,\delta + \eta(q)\eqno (3.47)$$ \vskip .5pc Let $\sigma\in
W_F^+$. Then by (3.20) one has $\lambda^{\sigma}\in X$, i.e. $D_{alcove}\s X$, and in addition (3.20)
defines a distinguished positive root partition of $\lambda^{\sigma}$. We denote this partition by
$q^{\sigma}$ and refer to this partition as the $\sigma$ positve root partition of
$\lambda^{\sigma}$. Thus $q^{\sigma}\in Q_{\lambda^{\sigma}}$ is given by $q^{\sigma}_i =
n_{\varphi_i}(\sigma)$ for all $i\in I_m$, using the notation of \S1.3. The following result
characterizes the elements in the subset $D_{alcove}\s X$ and for, each $\sigma\in W_f^+$, the
result characterizes the $\sigma$ positive root partition among all the positive root partitions in
$Q_{\lambda^{\sigma}}$. \vs {\bf Theorem 3.14.} {\it Let
$q\in Q$. Then
$$c(q) \geq (\mu + \rho,\mu + \rho) - (\rho,\rho)\eqno (3.48)$$ where $\mu = \eta(q)$. Furthermore
one has equality in (3.48) if and only if $\mu = \lambda^{\sigma}$ for some (necessarily
unique, see (2.8)) $\sigma\in W_f^+$ and
$q = q^{\sigma}$.} \vs {\bf Proof.} Let $k = c(q)$ so that $Z\s (\wedge\u^-)_k$ where $Z$ is
$\g_{\d}$-submodule generated by $z^{(q)}$. But then, using the notation of (3.29), one readily has,
using e.g. the Freudenthal result (5.9.2) in [Ko-1], 
$(\mu + \rho,\mu + \rho) - (\rho,\rho) \leq m_k$. But $m_k\leq k$ by (3.29). This establishes the
inequality (3.38). 

Now if $\mu = \lambda^{\sigma}$ and $q = q^{\sigma}$ one has equality in (3.48) by (3.40).
Conversely if one has equality in (3.48) then, by Theorem 3.9, $Z$ is a $\g_{\d}$-irreducible
component of $Har(\wedge\u^-)$ and $z^{(q)}$ is a highest weight vector. Then by Theorems 3.10 and
3.13 there exists $\sigma\in W_f^+$ such that the highest weight of $Z$ is
$-Cas(\lambda^{\sigma})\,\delta + \lambda^{\sigma}$ and this weight occurs with multiplicity 1 in
$\wedge\u^-$. But then $z^{(q)} = z^{\sigma}$ (up to scalar multiplication). But then $q =
q^{\sigma}$ by Remark 3.12. QED\vs 3.8. Using, in the present context, notation introduced in \S5.1
of [Ko-1], one defines an operation $\buildrel .\over +$, referred to as root addition, on the set
of all subsets of the set of affine roots. If $\Psi_i\s\www\Delta',\,i=1,2$, then $\Psi = \Psi_1
\buildrel .\over + \Psi_2$ if $\Psi = \{\beta\in\www\Delta'\mid\beta =
\beta_1 +\beta_2,\,\hbox{for some}\,\beta_i\in \Psi_i\}$. Let $\Psi\s \www\Delta'$. We will say that
$Psi$ is closed under root addition if $\Psi \buildrel .\over + \Psi = \Psi$ and $\Psi$ is abelian or
commutative if $\Psi \buildrel .\over + \Psi = \emptyset$. A subset $\Phi\s\Delta_+$ is called
ideal in $\Delta_+$ if $\Delta_+ \bd \Phi \s \Phi$. If $\Phi_i,\,i=1,2,$ are two such ideals then
obviously $\Phi_1\bd\Phi_2$ is again such an ideal. 

Let $ \sigma\in W_f^+$ and let $L(\sigma) = max_{\varphi\in \Delta_+}n_{\varphi}(\sigma)$ using
the notation of \S 1.3. Since $\psi$ is the highest weight of the adjoint representation, clearly
(see Remark 1.2) 
$$n_{\psi}(\sigma)\geq n_{\varphi}(\sigma)\eqno (3.49)$$ for any $\varphi \in \Delta_+$. Thus
$$L(\sigma) = n_{\psi}(\sigma) \eqno (3.50)$$ which, in the notation of \S1.4 implies that
$$A_{\sigma}\s \hh^{(L(\sigma) +1)}\,\,\hbox{but}\,\,A_{\sigma}\not\subset \hh^{(L(\sigma))}\eqno
(3.51)$$ In Theorem 3.16 below we observe that
$\sigma$ defines a chain of $L(\sigma)+1$ ideals (not necessarily distinct) of $\Delta_+$. For any
$i\in \Bbb Z_+$ let
$$\Delta_i(\sigma)=
\{\varphi\in
\Delta_+\mid i\leq n_{\varphi}(\sigma)\}$$ so that $\Delta_i(\sigma) = \emptyset $ if 
$i>L(\sigma)$. Clearly $$\Delta_{L(\sigma)}(\sigma)\s\cdots \s\Delta_0(\sigma)=\Delta_+\eqno (3.52)$$
\vskip .5pc {\bf Remark 3.15.} Observe that if $i\in \Bbb N,\,\varphi\in \Delta_+$ and
$\sigma\in W_f^+$ then, by Proposition 3.1,
$$\varphi \in \Delta_i(\sigma) \iff i\,\delta- \varphi\in \Phi_{\overline {\sigma}}\eqno
(3.53)$$ (Note the exclusion of $i=0$). \vskip 1pc If $\varphi_1,\varphi_2\in\Delta_+$ and
$\varphi_1 + \varphi_2$ is a root then note that by Remark 1.2 one has $$n_{\varphi_1 +
\varphi_2}(\sigma) \in \{n_{\varphi_1}(\sigma) + n_{\varphi_2}(\sigma), n_{\varphi_1}(\sigma) +
n_{\varphi_2}(\sigma) +1\}\eqno (3.54)$$ for any $\sigma\in W_f^+$. \vs {\bf Theorem 3.16.} {\it Let
$\sigma\in W_f^+$. Then
$\Delta_i(\sigma)$ is an ideal in $\Delta_+$ for any $i\in \Bbb Z_+$. Furthermore
$$\Delta_i(\sigma)\bd
\Delta_j(\sigma) \s
\Delta_{i + j}(\sigma)\eqno (3.55)$$ for any $i,j\in \Bbb Z_+$ so that in particular
$\Delta_{L(\sigma)}(\sigma)$ is an abelian ideal in $\Delta_+$. Finally using the
notation of \S 3.2 one has $$\lambda^{\sigma} = \sum_{i\in I_{L(\sigma)}}\langle
\Delta_i(\sigma)\rangle\eqno (3.56)$$} \vs {\bf Proof.} The first statement and (3.55) are
 immediate consequences of (3.54). On the other hand (3.56) clearly follows from (3.20) and
(3.53). QED \vskip1.5pc
\centerline{\bf 4. The main results}\vskip 1pc 4.1. We recall results in the 1965 paper
[Ko-2].  Let
$\u\s
\g$ be  any (complex) subspace. Let $k = \hbox{dim}\,\u$ so that $\wedge^k\u$ is a 1-dimensional
subspace of
$\wedge^k\g$. Let $M$ be the maximal dimension of an abelian subalgebra
of $\g$. For any $k\in \Bbb Z_+$ let $C_{k}\s \wedge^k\g$ be the span
of all 1-dimensional subspaces of the form $\wedge^k \a$ where $\a\s
\g$ is a $k$-dimensional abelian subalgebra of 
$\g$. Obviously $C_k = 0$ if $k>M$. Clearly $C_{k}$ is a $\g$-submodule of 
$\wedge^k\g$ under the adjoint action of $\g$. Of course
the same is true of
$$C = \sum_{\k=0}^M C_k$$ 

If a (complex) subspace $\u\s\g$ is stable under 
$\hh$ 
 let $\Delta(\u) = \{\varphi\in \Delta\mid e_{\varphi}\in
\u\}$. Let $\b$ be the Borel subalgebra of $\g$ containing $\hh$ such that $\Delta(\b)=
\Delta_+$ and let $\n = [\b,\b]$ be the nilradical of $\b$. Let $\Xi$ be an index set
parameterizing the set of all abelian ideals in $\b$ and for any
$\xi\in \Xi$ let $\a_{\xi}$ be the corresponding abelian ideal. For any $\xi\in \Xi$ it is
immediate that $\a_{\xi}\s\n$ so that $\Xi$ is finite and 
$$\a_{\xi} = \sum_{\varphi\in \Delta(\a_{\xi})} \Bbb C e_{\varphi}\eqno
(4.1)$$ and hence $$\wedge^k\a_{\xi} = \Bbb
C\,e_{\varphi_1}\wedge\cdots\wedge e_{\varphi_k}\eqno (4.2)$$ where $k =
\hbox{dim}\,\a_{\xi}$ and $$\Delta(\a_{\xi}) =
\{\varphi_1,\ldots,\varphi_k\}\eqno (4.3)$$ The subsets of $\Delta_+$ of the form (4.3) are
characterized by Theorem (8) in [Ko-2]. For $k\in \Bbb Z_+$ where $k\leq card\,\Delta_+$ let
$\Xi_k= \{\xi\in \Xi\mid dim\,\a_{\xi} = k\}$. Of course $\Xi_k$ is empty if $k>M$. If $\mu\in \hh$
let
$|\mu| = (\mu,\mu)^{1\over 2}$. The characterization is as follows:
\vs {\bf Theorem 4.1.} {\it Let $\Phi\s\Delta_+$ and let $k= card\,\Phi$. Then if $\Phi =
\{\varphi_1,\ldots,\varphi_k\}$ one has $$|\rho + \varphi_1+\cdots+ \varphi_k|^2- |\rho|^2\leq
k\eqno (4.4)$$ and one has equality in (4.4) if and only if there exists $\xi\in \Xi_k$ (necessariy
unique) such that $\Phi = \Delta(\a_{\xi})$.}\vs The adjoint action of $\g$ on itself induces
the structure of a $\g$-module (and hence $U(\g)$-module) on $\wedge\g$. Let $\xi\in \Xi$.
Since
$[\b,\a_{\xi}]\s \a_{\xi}$ it is immediate from (4.2) that, if $\xi\in \Xi_k$, then $\wedge^k
\a_{\xi}\s
\wedge^k\g$ is a highest weight space and hence, under the action of $\g$, generates an
irreducible
$\g$-submodule, denoted here by $V_{\xi}$, of $\wedge^k\g$.
Furthermore if we put $$\lambda_{\xi} = \sum_{\varphi\in \Delta(\a_{\xi})}
\varphi\eqno (4.5)$$ then $\lambda_{\xi}$ is the highest weight of
$V_{\xi}$. Theorem (7) in [Ko-2] implies that, for $\xi,\xi'\in \Xi$, $$\lambda_{\xi}
=\lambda_{\xi'}\iff \xi= \xi'\eqno (4.6)$$ This accounts for the multiplicity-free statement 
in the following result (Theorem 4.2). Like Theorem 4.1, Theorem 4.2 is also part of Theorem (8) in
[Ko-2].
\vskip 1pc {\bf Theorem 4.2.} {\it For any $\k\in \Bbb Z_+$ where $0\leq k\leq
M$ one has the direct sums 
$$C_k = \sum_{\xi\in\, \Xi_k} V_{\xi}\eqno (4.7)$$ and $$C = \sum_{\xi\in
\,\Xi} V_{\xi}$$ Furthermore $C$ and, a fortiori, $C_k$, are
multiplicity-free $\g$-modules.} \vs 4.2.  It is a beautiful result of Dale Peterson that
$card\,\Xi = 2^{\ell}$. Although Peterson's proof has not been published a sketch of his proof
appears in \S 2 of [Ko-5]. A key (and, for me, surprising) point of the proof was the
connection established between
$\Xi$ and a subset of $W_f^+$. Expanding on this connection P. Cellini and P. Papi
published a  simpler proof of Peterson's theorem in [C-P]. See Theorem 2.9 in that reference. 
Recalling the notation of 
\S 1.4 in our present paper here let 
$W_f^{(k)} = \{\sigma\in W_f^+\mid A_{\sigma}\s \hh^{(k)}\}$. By Proposition 1.4 one has
$card\,W_f^{(2)} = 2^{\ell}$. Peterson's theorem follows from a bijection $$\Xi\to W_f^{(2)}\eqno
(4.8)$$ established in [C-P]. \vs {\bf Remark 4.3.} Note that if $\sigma\in W_f^+$ then $$\sigma\in
W_f^{(2)}\iff n_{\varphi}(\sigma)\in \{0,1\},\,\, \forall \varphi\in \Delta_+\eqno (4.9)$$ Indeed by
definition $\sigma\in W_f^{(2)}$ if and only if $n_{\psi}(\sigma)\in \{0,1\}$. But then (4.9)
follows from (3.49). \vs Recently Ruedi Suter in [Su] showed (again a surprise for me) showed that
Peterson's theorem, in fact, follows from a result in [Ko-2].  In more detail, Theorem 4.4, below,
was known to me before [Su]. In fact it is an immediate consequence of Peterson's Proposition 2.5 in
[Ko-5] and the result (4.8) in [C-P]. The equation (4.10) below was discovered independently by Suter
in [Su]. But the main novelty is that the proof of (4.10),(4.11) and hence (4.13) in [Su] 
depends only on a 1965 result in [Ko-2], stated in the present paper as Theorem 4.1. With the benefit
of this knowledge we will prove Theorem 4.4 using only [Ko-2] and results established in the present
paper. 
\vs {\bf Theorem 4.4.} {\it For any
$\xi\in
\Xi$ there exists an (necessarily unique) element $\sigma_{\xi}\in W_f^+$ such that $$\lambda_{\xi} 
= \lambda^{\sigma_{\xi}}\eqno (4.10)$$ Moreover $\sigma_{\xi}\in W_f^{(2)}$ and the map $$\xi\to
W_f^{(2)},\qquad \xi\mapsto \sigma_{\xi}\eqno (4.11)$$ is a bijection. In particular one has
the inclusion $$\{\lambda_{\xi}\mid \xi\in \Xi\} \s D_{alcove}\eqno (4.12)$$ and the count
(Peterson's theorem) $$card\,\Xi = 2^{\ell}\eqno (4.13)$$}\vs {\bf Proof.} Let $\xi\in \Xi$. Then,
in the notation of \S3.7, (4.5) is a root partition $q$ of $\lambda$ (see (3.43)) and
$q_i\in\{0,1\}$ for any $i\in I_{m}$. But then if $\xi\in \Xi_k$ one has $c(q) = k$ (see (3.45)).
But then $c(q) = |\lambda_{\xi} + \rho|^2 - \rho|^2$ by Theorem 4.1. Thus there exists an
element $\sigma_{\xi}\in W_f^+$ satisfying (4.10) and $q= q^{\sigma_{\xi}}$ by Theorem 3.14.
Obviously if $\varphi\in \Delta_+$ then $$n_{\varphi}(\sigma_{\xi}) = 1\,\,\hbox{if}\,\,\varphi\in
\Delta(\a_{\xi})\,\,\hbox{and}\,\,n_{\varphi}(\sigma_{\xi}) = 0\,\,\hbox{if}\,\,\varphi\notin
\Delta(\a_{\xi})\eqno (4.14)$$ But then $\sigma_{\xi}\in W_f^{(2)}$ by (4.9). The map (4.11) is
injective by (4.6). 

Conversely let $\sigma\in W_f^{(2)}$. Without loss we may assume that $\sigma\neq 1$. Then
by (4.9) and Theorem 3.16 the set $$\Phi = \{\varphi\in \Delta_+\mid n_{\varphi}(\sigma)=1\}\eqno
(4.15)$$ is an abelian ideal in $\Delta_+$. Thus there exists $\xi\in \Xi$ such that $\lambda_{\xi} =
\langle\Phi\rangle$. On the other hand $\lambda^{\sigma} = \langle\Phi\rangle$ by (3.20). Hence
$\sigma = \sigma_{\xi}$. Consequently (4.11) is surjective. QED\vs As a consequence of Remark 4.3
and Theorem 4.4 one establishes the following property of $W_f^{(2)}$. \vs {\bf Theorem 4.5.} {\it
Let $\sigma\in W_f^+$. Then $$Cas(\lambda^{\sigma}) \geq \ell(\sigma)\eqno (4.16)$$ and equality
occurs in (4.16) if and only if $\sigma\in W_f^{(2)}$. Furthermore in that case writing $\sigma =
\sigma_{\xi}$ for $\xi\in \Xi$ (Theorem 4.4) one has $$\eqalign{Cas(\lambda^{\sigma}) &=
\ell(\sigma)\cr &= dim\,\a_{\xi}\cr}\eqno (4.17)$$}\vs {\bf Proof.} The inequality (4.16) follows
from (1.6) and (3.40). But then equality occurs if and only if $n_{\varphi}(\sigma)\in \{0,1\}$ for
all $\varphi\in \Delta_+$. But this is the case if and only if $\sigma\in W_f^{(2)}$ by Remark
4.3. The final statement then follows from Theorem 4.4 and (4.14). QED \vs 4.3. Recall 
\S 3.3. Define a pairing of $\u^-$  and $\u = t\,\g[t]$ so that for $p,q\in \Bbb N$ and
$x,y\in\g$ then
$(t^{-p}\,x,t^{q}\,y) = 0 $ if $p\neq q$ and $(t^{-p}\,x,t^{p}\,y) = (x,y)$. Let $y_k,k\in J$ be
the basis of $\g$, dual to the $x_j$, so that $y_k$ is a weight vector of weight $-\mu_k$. It
follows then that $\{t^q\,y_k \mid (q,k)\in \Bbb N\times J\}$ is an $\www\hh$ basis of $\u$, dual
to the basis $\{t^{-p}\,x_j \mid (p,j)\in \Bbb N\times J\}$ of $\u^{-}$. We may identify $\u^-$
here with the subspace of all linear functionals $f$ on $\u$ which vanish on $t^{N}\g$ (where $N$
depends on
$f$) for sufficiently large $N$.
 
If $r$ is defined as in (3.23) let $w_r\in\wedge^n\u$ be defined so that $$w_r =
t^{r_1}y_{r_{[1]}}\wedge \cdots \wedge t^{r_n}y_{r_{[n]}}\eqno (4.18)$$ The obvious analogue of
Proposition 3.5 and 3.6 is \vs {\bf Proposition 4.6.} {\it Let $r\in R_n$. Then $w_r$ is a
$\www\hh$-weight vector of weight $$p(r)\delta-\mu(r)$$ and the set $\{w_r\mid r\in R_n\}$ is a
basis of $\wedge^n\u$.} \vs The pairing of $\u^-$ and $\u$ extends, as usual (determinently), to a
nonsingular pairing of $\wedge\u^-$ and $\wedge \u$. The subspaces $\wedge^m\u^-$ and $
\wedge^n\u$ are orthogonal if $m\neq n $ and if $n=m$ one notes that $\{w_r\mid r\in R_n\}$ and
$\{z_r\mid r\in R_n\}$ are dual bases. The pairing of $\wedge\u^-$ and $\wedge \u$ is clearly
invariant under the action of $\g_{\d}$. If $y\in \u$ let $\theta(y)\in End\,\wedge\u$ be the
operator of the adjoint action of
$y$ on $\wedge\u$ and if $x\in \u^-$ let $\iota(x)\in End\,\wedge\u$ be the interior product by $x$.
Let $\partial_+$ be the boundary operator on $\wedge \u$ whose derived homology is $H_*(\u)$. Using
a standard expression for $\partial_+$ one has that if $u\in \wedge\u$ then $$\partial_+\,u =
1/2\sum_{(p,j)\in \Bbb N\times J}\theta(t^p\,y_j)\,\iota(t^{-p}\,x_j)\,u\eqno (4.19)$$ noting that
$\iota(t^{-p}\,x_j)\,u$ is nonzero for only a finite subset of $\Bbb N\times J$. Let 
$\theta^*(t^p\,y_j)\in End\,(\wedge\u)^*$ be the negative transpose of $\theta(t^p\,y_j)$. It is
immediate that $\wedge\u^-$ is stable under $\theta^*(t^p\,y_j)$ and that, using the notation of
(3.26), $$\theta^*(t^p\,y_j):(\wedge\u^-)_k\to (\wedge\u^-)_{k-p}\eqno (4.20)$$ noting
$(\wedge\u^-)_j= 0$ if $j$ is negative. For any $x\in \u^-$ let $\varepsilon(x)\in
End\,\wedge\,\u^-$ be the operator of exterior multiplication by $x$. One sees that if
$d_+\in End\,(\wedge\u)^*$ is the negative transpose of $\partial_+$ then $\wedge\u^-$ is stable
under $d_+$ and, in fact for any $v\in 
\wedge\u^-$ one has $$d_+\,v = 1/2 \sum_{(p,j)\in \Bbb N\times
J}\varepsilon(t^{-p}\,x_j)\theta^*(t^p\,y_j)\,v\eqno (4.21)$$ noting that, by (4.20),
$\theta^*(t^p\,y_j)\,v\neq 0$ for only a finite subset of $\Bbb N\times J$. Of course the pair
$((\wedge\u)^*,d_+)$ is the cochain complex whose derived cohomology is $H^*(\u)$. The pair
$(\wedge\u^-,d_+)$ is a subcomplex and we denote the derived cohomology by ${\cal H}^*(\u)$. Since
$d_+$ is antiderivation in either case both $H^*(\u)$ and ${\cal H}^*(\u)$ are algebras and one has
an algebra homomorphism
$${\cal H}^*(\u)\to H^*(\u)$$ On the other hand (see \S 3.4) $\wedge\u^-$ has a bigrading
$(\wedge^n\u^-)_k$ and one notes from (4.21) that $$d_+:(\wedge^n\u^-)_k\to
(\wedge^{n+1}\u^-)_k\eqno (4.22)$$ In particular $((\wedge\u^-)_k,d_+)$  (see
\S 3.4.) is a finite dimensional $\g_{\d}$-completely reducible subcomplex of $(\wedge\u^-,d_+)$, for
$k\in \Bbb Z_+$. The derived cohomology is denoted by ${\cal H}(\u)_k$. Since the $r_i$ in (3.25)
are positive one notes that $(\wedge^n\u^-)_k= 0$ for
$n>k$. Also note that (see (3.25)) for $n,n',k,k'\in \Bbb Z_+$, $$(\wedge^n\u^-)_k\wedge
(\wedge^{n'}\u^-)_{k'}\s (\wedge^{n+n'}\u^-)_{k+k'}\eqno (4.23)$$ Of course $\wedge$ induces cup
product ($\vee$) in ${\cal H}(\u)$. This establishes
\vs {\bf Proposition 4.7.} {\it One has the direct sum $${\cal H}(\u)= \sum_{k\in \Bbb Z_+}
({\cal H}(\u))_k\eqno (4.24)$$ where $({\cal H}(\u))_k$ is a finite dimensional $\g_{\d}$-completely
reducible $\g_{\d}$-module. Furthermore with regard to cohomological degree $$({\cal
H}(\u))_k =
\sum_{n=0}^k({\cal H}^n(\u))_k\eqno (4.25)$$ and one has the cup product relation
 $$({\cal H}^n(\u))_k\vee 
({\cal H}^{n'}(\u))_{k'}\s ({\cal H}^{n+n'}(\u))_{k+k'}\eqno (4.26)$$}

\vs

\def\wd{\wedge\u^-}

\rm 
4.4. In \S3.4 we introduced Garland's harmonic subspace
$Har(\wd)$ of $\wedge\u^-$. The subspace $Har(\wd)$ is a space of
cycles for the boundary operator on $\wd$ and the quotient map
$Har(\wd) \to H_*(\u^-)$ is an isomorphism (see (3.28)). We will
now see that $Har(\u^-)$ also represents the cohomology ${\cal
H}(\u)$.  We will first clarify the structure of $Har(\wd)$. For
any $\sigma\in W_f^+$ we have defined a $\g_{\d}$-irreducible
submodule $Z_{\sigma}$ of $\wd$ in (Garland) Theorem 3.10.
Combining Theorem 3.10 and Theorem 3.13 one has \vs {\bf Theorem
4.8.} {\it The highest weight of $Z_{\sigma}$ is
$-Cas(\lambda^{\sigma})\delta + \lambda^{\sigma}$. Furthermore this
weight occurs with multiplicity 1 in $\wd$ so that, a fortiori, the
representation of $\g_{\d}$ afforded by $Z_{\sigma}$ occurs with
multiplicity 1 in $\wd$. Next
$$Z_{\sigma}\s
(\wedge^{\ell(\sigma)}(\u^-))_{Cas\,\lambda^{\sigma}}\eqno (4.27)$$
In fact $$(Har_n(\wd))_k = \sum_{\sigma\in
W_f^+,\,\ell(\sigma)=n,\,Cas(\lambda^{\sigma})=k} Z_{\sigma}\eqno
(4.28)$$ }\vs If $\lambda\in D$ then, as one knows, the dual
$(V_{\lambda})^*$, as a $\g$-module, is characterized by the property that $-\lambda$
is the lowest weight of
$(V_{\lambda})^*$. A similar statement is clearly true for the 
reductive Lie algebra $\g_{\d}$. As a consequence of Theorem 4.8 one can then make
the following \vs {\bf Remark 4.9.} If $\sigma\in W_f^+$ then the dual
 $(Z_{\sigma})^*$, as a $\g_{\d}$-module, is characterized by the property that
$Cas(\lambda^{\sigma})\delta - \lambda^{\sigma}$ is the lowest weight of
$(Z_{\sigma})^*$. \vs  Now one knows that there exists an automorphism
$\theta$ on
$\g$ which stabilizes $\hh$ and is such that $\theta|\hh$ is minus the identity. In
particular
$\theta(\Delta_+) = -\Delta$ so that if $\kappa\in W$ is the long element then
$\kappa\theta$ stabilizes $\Delta_+$ and also stabilizes $D$. If $\lambda\in D$ let
$\lambda'= \kappa\theta(\lambda)$. Since the $-\lambda$ is the lowest weight of
$V_{\lambda'}$ we may identify $V_{\lambda'}$ with the dual $\g$-module
$V_{\lambda}^*$. \vs {\bf Remark 4.10.} As an application of Theorem 2.4 note that
$D_{alcove}$ is stable under $\kappa\theta$. Indeed since
$\chi_{\lambda}(a_P)\in \{-1,1\}$ clearly (by the reality of the character
value)
$\chi_{\lambda'}(a_P)= \chi_{\lambda}(a_P)$ so that $\chi_{\lambda'}(a_P)\in
\{-1,1\}$. Thus there exists an involutory bijection $W_f^+\to W_f^+,\quad
\sigma\mapsto \sigma'$ such that $$(\lambda^{\sigma})' = \lambda^{\sigma'}\eqno
(4.29)$$ Examples exist where $\sigma'\neq \sigma$.\vs The automorphism $\theta$
clearly extends to an automorphism of $\www\g$ which stabilizes $\www\hh$ and is
minus the identity on $\www\hh$. But then $\theta$ interchanges $\u$ and $\u^-$ since
$\theta(\delta) = -\delta$. As noted in Remark 4.9, for $\sigma\in W_f^*$, the
irreducible
$\g_{\d}$-module with lowest weight $Cas(\lambda^{\sigma})\delta -\lambda^{\sigma}$
readily identifies with $(Z_{\sigma})^*$, the dual
$\g_{\d}$-module to $Z_{\sigma}$. But now
Theorem 3.10 and Theorem 3.12 determines the $\g_{\d}$-module structure of $H_*(\u)$.
Recalling the notation at the end of \S 3.3 one has the direct sum  $$H_*(\u)
=\sum_{k\in \Bbb Z_+} (H_*(\u))_{-k}\eqno (4.30)$$ Furthermore $(H_*(\u))_{-k}$ is a
finite dimensional completely reducible $\g_{\d}$-module and as such $$(H_*(\u))_{-k}
\equiv \sum_{\sigma\in W_f^+,\,Cas(\lambda^{\sigma})= k} (Z_{\sigma})^*\eqno
(4.31)$$\vs We can now prove that $Har(\wd)\s \wd$ represents the cohomology 
${\cal H}(\u)$ as well as the homology  (see (3.28)) $H_*(\u^-)$. \vs {\bf Theorem
4.11.} {\it Any element in $Har(\wd)$ is a $d_+$-cocycle and the induced linear map
$$Har(\wd) \to {\cal H}^*(\u)\eqno (4.32)$$ is a $\g_{\d}$-module isomorphism. In
particular (4.32) restricts to an isomorphism $$(Har(\wd))_k \to ({\cal
H}^*(\u))_k\eqno (4.33)$$ of finite dimensional completely reducible
$\g_{\d}$-modules, for any $k\in
\Bbb Z_+$.}\vs {\bf Proof.} It clearly suffices to prove (4.33). But the
nonsingular pairing of $\wd$ and $\wedge\u$ induces a nonsingular pairing of the finite
dimensional completely reducible $\g_{\d}$-modules $(\wd)_k$ and $(\wedge\u)_{-k}$.
But $d_+|(\wd)_k$ is the negative transpose of $\partial_+|(\wedge\u)_{-k}$. Thus, as
$\g_{\d}$-modules one has $$({\cal H}^*(\u))_k \equiv \sum_{\sigma\in
W_f^+,\,Cas(\lambda^{\sigma})= k} Z_{\sigma}$$ by (4.31). But by complete
reducibility and the multiplicity 1 statement in Theorem 4.8 one must have
$Har_n(\wd))_k\s Ker d_+$ (see (4.28)) and the $\g_{\d}$-isomorphism (4.33). QED \vs
4.5. We now introduce a new grading ${\cal H}^{[j]}(\u)$ on ${\cal H}^*(\u)$. In
doing so we are following a suggestion of Pavel Etingof who pointed out to us that
our subsequent results can be neatly formulated using this grading. As noted in \S
4.3 one has $(\wedge^n\u^-)_k = 0$ for $n>k$. In particular $({\cal H}^n(\u))_k = 0$
for $n>k$ (see (4.25)). For $j\in \Bbb Z$ let $$\wedge^{[j]}\u^- = \sum_{n,k\in
\Bbb Z_+,\,k-n = j}(\wedge^n\u)_k\eqno (4.34)$$ and let $${\cal H}^{[j]}(\u) =
\sum_{n,k\in
\Bbb Z_+,\,k-n = j}({\cal H}^n(\u))_k\eqno (4.35)$$ so that one has the direct
sums 
$$\eqalign{\wd &= \sum_{j\in \Bbb Z_+}\wedge^{[j]}\u\cr
{\cal H}^*(\u) &= \sum_{j\in \Bbb Z_+}{\cal
H}^{[j]}(\u)\cr}\eqno (4.36)$$
\vskip .5pc 
{\bf Proposition 4.12.}
{\it The subspaces $\wedge^{[j]}\u^-$ and, a fortiori, ${\cal H}^{[j]}(\u)$ are finite
dimensional. Moreover, with respect to wedge and cup product
$$\eqalign{(\wedge^{[j]}\u^-)\wedge (\wedge^{[j']}\u^-)&\s (\wedge^{[j+j']}\u^-)\cr
{\cal H}^{[j]}(\u)\vee {\cal
H}^{[j']}(\u)&\s {\cal H}^{[j + j']}(\u)\cr}\eqno (4.37)$$ 
In particular $\wedge^{[0]}\u^-$ is a finite dimensional subalgebra of $\wd$ and
${\cal H}^{(0)}(\u)$ is a finite dimensional subalgebra of
${\cal H}^*(\u)$. } \vs {\bf Proof.} Clearly (4.37) follows from (4.23). It
suffices only to show that $(\wedge^{[j]}\u^-)$ is finite dimensional for any
$j\in \Bbb Z_+$. Recalling the notation of
\S3.3 note that, for any $r\in R_n$, $$z_r \in \wedge^{[|p(r)| - n]}\u^-\eqno (4.38)$$
But, by (3.24), the number of
$i\in I_n$ such that
$r_i = 1$ is at most
$dim\,\k$ for any $n$. This implies that $|p(r)| - n \geq n- dim\,\k$. But this readily
implies that $\wedge^{[j]}\u^-$ is finite dimensional, for any $j\in \Bbb Z_+$.
QED\vs {\bf Remark 4.13}. Note that as a $\g_{\d}$-module one has $${\cal
H}^{[j]}(\u)\equiv \sum_{\sigma\in W_f^+, Cas(\lambda^{\sigma})- \ell(\sigma)=
j}Z_{\sigma}\eqno (4.39)$$ Indeed (4.39) follows from (3.33), (3.35), (3.36) and
(3.39). \vs  Now note that by (4.22) one has
$$d_+:\wedge^{[j]}\u^-\to
\wedge^{[j-1]}\u^-\eqno (4.40)$$ for any $j\in \Bbb Z$. This implies part of \vs {\bf
Lemma 4.14.} {\it One has $$\wedge^{[0]}\u^- \s Ker \,d_+\eqno (4.41)$$ Moreover as an
algebra and a
$\g_{\d}$ module one has 
$\wedge^{[0]}\u^-\equiv \wedge\g$. In fact there exists a $\g_{\d}$-module algebra
isomorphism $$\wedge\g \to \wedge^{[0]}\u^-\eqno (4.42)$$ where, if
$u_i\in\g,\,i=1,\ldots,k,$ then $$u_1\wedge\cdots\wedge u_k \mapsto
t^{-1}u_1\wedge\cdots\wedge t^{-1}u_k\eqno (4.43)$$}\vs {\bf Proof.} One has (4.41)
by (4.40) since of course $\wedge^{[-1]}\u^- = 0$. Recalling (3.25) the condition that
$z_r\in \wedge^{[0]}\u^-$ is that $r_i =1 $ for all $i\in I_n$. But this clearly
implies that there is a $\g_{\d}$-module algebra isomorphism (4.42) satisfying
(4.43). QED \vs Now if $u\in \g$ then in the notation of (4.21)), clearly
$\theta^*(t^p y_j) (t^{-2}u)= 0$ if $p\geq 2$ (since $\wd$ is stable, by definition,
under $\theta^*(t^p y_j)$) and $$\varepsilon(t^{-1}x_j)\theta^*(t^1 y_j) (t^{-2}u)=
t^{-1}x_j\wedge t^{-1}[y_j,u]$$  Hence $$d_+(t^{-2}u)= \sum_{j\in J}t^{-1}x_j\wedge
t^{-1}[y_j,u]\eqno (4.44)$$ We have identified $\g$ with its dual $\g^*$ using the Killing
form so that if $d$ is the Cartan-Eilenberg-Koszul coboundary operator on $\wedge\g^*$ it, with
this identification, is an antiderivation of degree 1 in $\wedge\g$. Using the
standard formula for $d$ one then has $d\,u\in \wedge^2\g$ for any $u\in \g$ and
explicitly $$d\,u = 1/2\sum_{j\in J} x_j\wedge [y_j,u] \eqno (4.45)$$ Let $(d\,\g)$ be
the ideal in $\wedge \g$ generated by $d\,\g$. In \S4.1 we introduced the
multiplicity-free $\g$-submodule $C\s \wedge \g$ defined by the set of all abelian
subalgebras of $\g$. Theorem 4.2 asserts that the highest weight vectors in $C$ are
given by the $2^{\ell}$ abelian ideals in $\b$. Theorem 4.3. in [Ko-6] contains
the following statement \vs {\bf Theorem 4.15.} {\it One has the direct sum $$C \oplus
(d\g)\eqno (4.46)$$} \vs 
The space $C$ inherits an algebra structure, as a consequence of Theorem 4.15, since
(4.46) implies that $$C\equiv
\wedge\g/(d\,\g)\eqno (4.47)$$ At the time [Ko-6] was written I had no idea
about the meaning of this algebra structure. This question is resolved in the
following theorem (Theorem 4.16). We could use Theorem 4.15 to prove much of Theorem
4.16. However, we will, instead, prove Theorem 4.16 using
results established in the present paper. 

The advantage in dealing with cohomology ${\cal H}^*(\u)$ instead of homology
$H_*(\u^-)$ is that ${\cal H}^*(\u)$ has an algebra stucture. The
nature of this algebra is, nevertheless, presently, quite mysterious to us. However for
the finite dimensional subalgebra ${\cal H}^{[0]}(\u)$ one has the following result.
\vs {\bf Theorem 4.16.} {\it As a $\g$-module ${\cal H}^{[0]}(\u)$ is
multiplicity-free with $2^{\ell}$ irreducible components. In fact (recalling \S 4.1) 
$${\cal H}^{[0]}(\u)
\equiv \sum_{\xi\in \Xi} V_{\xi}\eqno (4.48)$$ As an algebra (under cup product)
$${\cal H}^{[0]}(\u)\equiv \wedge\g/ (d\g)\eqno (4.49)$$}\vs {\bf Proof.} The statement
(4.48) follows from (4.39) and Theorem 4.5. 

Now by (4.40) and (4.41) one has $${\cal H}^{[0]} =
\wedge^{[0]}\u^-/d_+(\wedge^{[1]}\u^-)\eqno (4.50)$$ But $z_r\in \wedge^{[1]}\u^-$, by
(3.25), if and only if all but one $r_i=1$ and the remaining $r_i$ equals 2. Thus
$$\wedge^{[1]}\u^- = (\wedge^{[0]}\u^-)\wedge t^{-2}\,\g\eqno (4.51)$$ But then
$d_+(\wedge^{[1]}\u^-)$ is the ideal in $\wedge^{[0]}\u^-$ generated by
$d_+(t^{-2}\,\g)$, by (4.41). The algebra isomorphism (4.49) then follows immediately
from (4.42), (4.43), (4.44) and (4.45). QED\vs Recall the notation of \S 4.1. If
$k\in \Bbb Z_+$ then we have defined $C_k\s \wedge^k\g$ in terms of all the abelian
subalgebras of $\g$ having dimension $k$. Theorem (5) in [Ko-2] gives a different 
characterization of $C_k$. Note that, by (2.1.7) in [Ko-2], the Laplacian $L$ is
$Cas/2$ operating in $\wedge\g$. For any $k\in \Bbb Z_+$ let $m_{(k)}$ be the
maximal eigenvalue of $Cas$ in $\wedge^k\g$. Theorem (5) in [Ko-2] then asserts
\vs {\bf Theorem 4.17.} {\it  For any $k\in \Bbb Z_+$ one has $$m_{(k)}\leq k\eqno
(4.52)$$ Furthermore equality occurs in (4.52) if and only if $k\leq M$ in which case
$C_k$ is the eigenspace for $Cas$ (operating in $\wedge^k\g$) corresponding to the
eigenvalue $k$.}
\vs {\bf Remark 4.18.}
 In the light of Lemma 4.14 and (4.48) it is not difficult to show that the 1965
result, Theorem 4.17 above, is, in fact, implied by Theorem 2.5 in Garland's 1975
paper [Gar].\vs 4.6. Returning to (2.9) let $b_k\in \Bbb Z$, for $k\in \Bbb Z_+$, be
defined so that $$(\prod_{n=1}^{\infty} (1-x^n))^{dim\,K} = \sum_{k\in \Bbb
Z_+}b_k\,x^{k}\eqno (4.53)$$ By (2.9) and (2.10) one has $$\sum_{k\in \Bbb
Z_+}b_k\,x^{k} = \sum_{\sigma\in
W_f^+}(-1)^{\ell(\sigma)}\,dim\,V_{\lambda^{\sigma}}\,x^{Cas(\lambda^{\sigma})}\eqno
(4.54)$$ which immediately yields the finite sum (see (4.16))  $$\eqalign{b_k &= 
\sum_{\sigma\in W_f^+, Cas(\lambda^{\sigma}) = k} (-1)^{\ell(\sigma)}
dim\,V_{\lambda^{\sigma}}\cr &= \sum_{\sigma\in W_f^+, Cas(\lambda^{\sigma}) = k}
 (-1)^{\ell(\sigma)}
dim\,Z_{\sigma}\cr}\eqno (4.55)$$ Note that the second equality in (4.55) follows
from the first line in Theorem 4.8. 

 Let $\v\s \n$ be the span of $\{e_{\varphi}\mid (\psi, \varphi) > 0\}$. 
One knows that $\v$ is a
Heisenberg Lie algebra so that we can write $\hbox{dim}\,\v = 2m+1$ where 
$m\in\Bbb Z_+$. One has
$\Bbb C\,e_{\psi}= \hbox{cent}\,\v$. Here we are regarding the case where $m=0$ (i.e.
when
$\g \equiv Lie\, Sl(2,\Bbb C)$ as a Heisenberg Lie algebra. From the Heisenberg
structure of
$\v$ one knows that there exists a partition $\Delta(\v) - \{\psi\} = 
\Delta^1(\v)\cup \Delta^2(\v)$
where each of the two parts has $m$ roots which can be ordered so that if $\Delta^1(\v) =
\{\beta_1,\ldots,\beta_m\}$ and $\Delta^2(\v) = \{\gamma_1,\ldots,\gamma_m\}$ then for
$i=1,\ldots,m$, $$ \beta_ i + \gamma_i = \psi\eqno (4.56)$$ 

If $\Bbb \g$ is simply laced (A-D-E case) then we have known for some time
that $m = h-2$ where
$h$ is the Coxeter number. See e.g. (1.10.1), p. 214 in [Ko-4]. In the
non-simply laced case we had checked that $m\geq \ell -1$. But one could do better. D.
Peterson informed us that $m= h^{\vee}-2$ in general where $h^{\vee}$ is the dual
Coxeter number (see \S1.5 and the notational change in Proposition 1.8). Knowing this
one readily supplies an easy proof.
\vskip 1pc {\bf Proposition 4.19.} {\it If $m$ is defined as in (4.56) then $m =
h^{\vee}-2$ where
$h^{\vee}$ is the dual Coxeter number (see \S1.5).} \vskip 1pc {\bf Proof.}  If
$\varphi\in
\Delta_+$ then $(\varphi,\psi)\geq 0$ since $\psi$ is the highest root. But if 
$\varphi\in (\Delta(\v) -\{\psi\})$ then $(\varphi,\psi) = (\psi,\psi)/2$ since $\psi$
is a long root. Thus $(2\rho,\psi) - (\psi,\psi) = m\,(\psi,\psi)$. Hence
$(2\rho,\psi) + (\psi,\psi)= (m+2)(\psi,\psi)$. But $(2\rho,\psi) + (\psi,\psi) = 1$
by (1.18). But then $m = h^{\vee}-2$ by (1.17). QED \vskip 1pc 
We can now prove \vskip 1pc {\bf
Theorem 4.20.} {\it Let $\sigma$ lie in the complement of $W_f^{(2)}$ in $W_f^+$ (see
\S 4.2). Then
$$\ell(\sigma) \geq h^{\vee}\eqno (4.57)$$ where $h^{\vee}$ is the dual Coxeter number
(see
\S1.5)).}\vskip 1pc  {\bf Proof.} By assumption, if $r^{\sigma}\in
Reg\,A_{\sigma}$ (see \S 1.2), one has
$\psi(r^{\sigma}) >2$. But then, using (4.56) and Proposition 4.19, for any $i =
1,\ldots,h^{\vee} -2,$ one must have either $\beta_i(r^{\sigma}) >1$ or
$\gamma_i(r^{\sigma}) >1$ and possibly both inequalities. In any case the number of
$\varphi$-walls, where $\varphi\in
\Delta(\u) - \{\psi\}$, separating
$r^{\sigma}$ and
$2\rho$ is at least $h^{\vee}-2$. But since $\psi(r^{\sigma}) >2$, the walls 
$\hh_{\psi,1}$ and
$\hh_{\psi,2}$ also separate $r^{\sigma}$ and $2\rho$. This accounts for 
$h^{\vee}$ separating walls. This
proves (4.57) (see (1.6)). QED \vskip 1pc As a corollary one has \vskip 1pc {\bf
Theorem 4.21.} {\it Let $\sigma\in W_f^+$. If $Cas(\lambda^{\sigma})\leq h^{\vee}$,
where
$h^{\vee}$ is the dual Coxeter number (see \S 1.5), then 
$\sigma \in W_f^{(2)}$ (see \S 4.2)) and $$Cas(\lambda^{\sigma}) = \ell(\sigma)\eqno
(4.58)$$}
 \vskip 6pt {\bf
Proof.} If $Cas(\lambda^{\sigma})< h^{\vee}$ then $\ell(\sigma)< h^{\vee}$ by (4.16).
Hence
$\sigma\in W_f^{(2)}$ by Theorem 4.20. If $Cas(\lambda^{\sigma}) = h^{\vee}$ then
$\ell(\sigma)\leq h^{\vee}$ by (4.16). But if
$\ell(\sigma) < h^{\vee}$ one has $\sigma\in W_f^{(2)}$ by Theorem 4.20. On the
other hand if
$\ell(\sigma) = h^{\vee}$ then one has the equality $\ell(\sigma) =
Cas(\lambda^{\sigma})$. But this implies that $\sigma\in W_f^{(2)}$ by Theorem 4.5.
In any case $\sigma\in W_f^{(2)}$. But then one has (4.58) by Theorem 4.5. QED\vskip
1pc Let
$P(t)$ be the power series given by (1.15) and defined by Bott, so that if
$P(t)=\sum_{k=0}^{\infty}p_k t^k$ then $p_k$ is the $2k$ Betti number of the loop
group $\Omega(K)$. Another consequence of Theorem 4.20 is that we can count the number
of abelian ideals in $\b$ whose dimension is $k$ when $k<h^{\vee}$. Recall the notation of
\S4.1. \vskip 1pc {\bf Theorem 4.22.} {\it If $k< h^{\vee}$ then $$\hbox{card}\,\, \Xi_k =
p_k\eqno (4.59)$$ where $p_k$ is the
$2k$ Betti number of the loop group $\Omega(K)$ and is given by (1.15).}
 \vskip 1pc {\bf Proof.} By (1.14) and 
(1.15) $ p_k$ is the number of alcoves $A_{\sigma},\,\sigma\in W_f^+$, in $\hh^+$ such that $\ell(\sigma)
= k$. But if
$k<h^{\vee}$ then, by Theorem 4.20, the set $\{\sigma\in W_f^+\mid \ell(\sigma)<h^{\vee}\}$ is contained
in
$W_f^{(2)}$. But then (4.59) follows from Theorem 4.3 and (4.17). QED\vskip 1pc Our main results
concern $b_k$ (see (4.53)) when $k\leq h^{\vee}$. One major difficulty in using (4.55) to determine
$b_k$ is the cancelation in the sums of (4.55) due to the alternation signs. When $k\leq h^{\vee}$ this
alternation disappears. The grading in $\wedge\g$ induces a grading in the quotient algebra
$\wedge\g/(d\g)$. $$\wedge\g/(d\g) = \sum_{k=0}^M \wedge^k\g/(d\g)^k\eqno (4.60)$$ where 
$(d\g)^k = (d\g)\cap\wedge^k\g$. The following theorem is one our main results. \vs {\bf
Theorem 4.23.} {\it  Assume
$k\leq h^{\vee}$. Then the following seven numbers are all equal $$\halign{#\thinspace & # \hfill \qquad
\qquad    & #\hfill\cr 
{\rm[1]} & $(-1)^kb_k$ & (see (4.53)) \cr
{\rm[2]} & $dim\, C_k$ & (see \S 4.1) \cr
{\rm[3]} & $ \sum_{\xi \in \Xi_k}dim\, V_{\xi}$ & (see \S 4.2) \cr
{\rm[4]} & $dim\, \{v \in \wedge^k \frak{g}\, \vert\,\theta(Cas)\,\,v = k\, v\}$ & (see [Ko-2])\cr
{\rm[5]} & $dim\, \wedge^k\g/(d\g)^k$ & (see (4.47)) \cr
{\rm[6]} & $dim\, {\cal H}^k(\u)$& (see \S 4.3)\cr
{\rm[7]} & $dim\, H_k(\u^-)$& (see \S3.4)\cr}$$}\vs {\bf Proof.} By Theorem 4.21 one may
replace the upper sum in (4.55) by $$b_k = (-1)^k\sum_{\sigma\in W_f^{(2)}, Cas(\lambda^{\sigma})=k}
dim\,V_{\lambda^{\sigma}}\eqno (4.61)$$ But then recalling the definition of $\Xi_k$ in \S 4.1 and
$V_{\xi}$ in \S 4.2 it follows from Theorem 4.3 and (4.17) that $$b_k = (-1)^k\sum_{\xi\in
\Xi_k}dim\,V_{\xi}\eqno (4.62)$$ This implies the equality of [1] and [3]. But then [2] and [3] are
equal by (4.7). But [2] and [4] are equal by Theorem 4.17. One has the equality of [2] and [5] by
(4.47). But both ${\cal H}^k(\u)$ and $H_k(\u^-)$ are in bijective correspondence with $Har_k(u_-)$ by
Theorem 4.11 and (3.28). In particular one has a linear isomorphism ${\cal H}^k(\u)\to H_k(\u^-)$.
This it suffices to prove the equality of [6] and [5]. But $$Har_k(\u^-) = \sum_{\sigma\in
W_f^+,
\ell(\sigma)= k} Z_{\sigma}\eqno (4.63)$$ by (4.28). But if $\ell(\sigma) = k$ then $\sigma\in
W_f^{(2)}$ by Theorem 4.20 and hence 
$Cas(\lambda^{\sigma}) = k$ by Theorem 4.5. Thus $$\eqalign{{\cal H}^k(\u)& = ({\cal H}^k(\u)_k\cr &\s
{\cal H}^{[0]}(\u)\cr}\eqno (4.64)$$ But then one has the equality of [6] and [5] by (4.49). QED\vs
{\bf Example}. Consider the case when $K = SU(5)$. Then $dim\,\k = 24$ so that $b_k = \tau(k+1)$ where,
using the terminology and notation of \S 4.5, Chapter 7, in [Se], $n\mapsto \tau(n)$ is the Ramanujan
tau function. In this case $h^{\vee} = 5$ and choosing, say [2] in
Theorem 4.23, the first 5 nontrivial Ramanujan numbers (see p. 97 in [Se] and also [L]) yield the
equality $$\eqalign{dim\,C_1& = 24\cr dim\,C_2 &= 252\cr dim\,C_3&= 1472\cr dim\,C_4& = 4870\cr
dim\,C_5&= 6048\cr}$$ \vskip .5pc We thank P. Etingof for pointing out to us that Theorem 4.23 (and its
proof) yield the acyclity of the complex (4.65) in Theorem 4.24 when $k\leq h^{\vee}$. Let $\partial_-$
be the boundary operator in
$\wedge\u$ whose derived homology is
$H_*(\u^-)$. Then, as noted in \S 3.4, if $k\in \Bbb Z_+$, $((\wedge\u^-)_k,\partial_-)$ is a finite
dimensional subcomplex of
$((\wedge\u^-)_k,\partial_-)$. The subcomplex is described by the $\partial_-$-maps
$$(\wedge^k\u^-)_k
\longrightarrow (\wedge^{k-1}\u^-)_k\longrightarrow\cdots \longrightarrow (\wedge^{0}\u^-)_k
\longrightarrow 0\eqno (4.65)$$ (Here we are using the fact that $(\wedge^n\u^-)_k = 0$ if $n>k$.
The latter statement is immediate from (3.25) and Proposition 3.6)\vs {\bf Theorem 4.24.} {\it If
$k\leq h^{\vee}$ then the complex (4.65) is acyclic. That is, $(H_n(\u^-))_k = 0$ unless $n=k$ so that 
$(H_*(\u^-))_k = (H_k(\u^-))_k$. In fact $(H_*(\u^-))_k = H_k(\u^-)$ and hence $$dim\,(H_*(\u^-))_k =
(-1)^kb_k\eqno (4.66)$$}\vs {\bf Proof.} One has $(H_n(\u^-))_k = 0$ unless $n=k$ by (4.28) and Theorem
4.21. But (4.28) implies that $(Har_k(\wedge\u^-))_k$ is, as a $\g$-module, given by the direct
sum
$\sum_{\xi\in \Xi_k}V_{\xi}$ (using Theorems 4.5 and 4.28). The remaining statements follow from
Theorem 4.23. QED\vs 4.7. For any complex number $s$ one can define $s$ power of Euler product
$\prod_{n=1}^{\infty} (1-x^n)$ by taking the logarithm of the Euler product, multiplying by $s$ and
then exponeniating. It follows easily that $$(\prod_{n=1}^{\infty} (1-x^n))^s = \sum_{k=0}^{\infty}
f_k(s)\,x^k\eqno (4.67)$$ where $f_k(s)$ is a polynomial of degree $k$ defined as follows: Let
$\mu:\Bbb N\to \Bbb Q$ be defined by putting $\mu(m) = \sum_{d|m}\,1/d$. For $k,n\in \Bbb
N,\,\,n\leq k$, let
$$Q_{k,n}=
\{q\in \Bbb N^n\mid q = (m_1,\dots,m_n),\,\,\sum_{i=1}^n m_i = k\}$$ and using this notation let
$$q_{k,n} = \sum_{q\in Q_{k,n}} \mu(m_1)\cdots \mu(m_n)$$ Put $f_0 = 1$. If $k\in \Bbb N$ let
$f_k(s)$ be the polynomial of degree $k$ (with 0 constant term ) defined by putting $$f_k(s) =
\sum_{n=1}^k q_{k,n}\,(-s)^n/n!\eqno (4.68)$$ Of course this
is a very complicated expression for
$f_k(s)$. In the notation of Theorem 4.23 one has $$ b_k = f_k(dim\,\k)\eqno (4.69)$$ Clearly the
polynomial
$f_k(s)$ would be known if we knew its roots.\vs {\bf Remark 4.25.} According to Serre (see top of p. 98
in [Se]) it is a question raised by D. H. Lehmer as to whether 24 is ever a root of $f_k(s)$ for any
$k\in
\Bbb Z_+$. \vs It is easy to see that $f_1(s) = -s$ and that in fact $0$ is a root of $f_k(s)$ for any
$k\in \Bbb N$. We will determine $f_2(s),\, f_3(s)$ and $f_4(s)$ in a uniform way using Theorem 4.23.
We first observe \vs {\bf Proposition 4.26.} {\it For the missing roots (to be determined in Theorem
4.27 below) $r_4,r_3$ and $r_2$ one has $$\eqalign{f_4(s) &= 1/4!\,\,s(s-1)(s-3)(s-r_4)\cr -f_3(s) &=
1/3!\,\,s(s-1)(s-r_3)\cr f_2(s) &= 1/2!\,\, s(s-r_2)\cr}$$}\vs {\bf Proof.} Euler has determined the
right side of (4.67) when
$s=1$. The only nonzero coefficients on the right side of (4.67) are the coefficients of the
pentagonal powers
$x^{(3n^2-n)/2}$ where $n\in \Bbb Z$. Since 3 and 4 are not pentagonal numbers it follows that 1 must
be a root of $f_3(s)$ and $f_4(s)$. Now Jacobi has determined the right side of (4.67) when $s=3$.
Here the only nonzero coefficients on the right side of (4.67) are the coefficients of the triangular
powers 
$x^{n(n+1)/2}$ for $n\in \Bbb Z_+$. Since $4$ is not a triangular number, $3$ must be a root of
$f_4(s)$. This proves the proposition. QED\vs Malcev has determined $M$ (see \S 4.1) for all complex
simple Lie algebras (see \S 4.3 in [Ko-2]). There are only 3 cases where $M<h^{\vee}$, namely the
cases where $\g$ is of $A_1,A_2$ and $G_2$. The relevant information is in the following table.
$$\matrix{\underline {\g\,type}&\underline {M}&\underline
{h^{\vee}} &\underline{dim\,\k}\cr 
A_1 & 1& 2& 3\cr A_2 & 2 &3&8 \cr G_2
&3 &4 &14\cr}$$ But then, by Theorem 4.23, $C_4=0$ and hence $b_4 = 0$ if $\g$ is of type $G_2$.
Next $C_3 = 0$ and hence $b_3 = 0$ if $\g$ is of type $A_2$. Finally $C_2=0$ and hence $b_2=0$ if
$\g$ is of type $A_1$. For these three cases $f_{h^{\vee}}(dim\,\k) = 0$ by (4.69). Hence we have
proved \vs {\bf Theorem 4.27.} {\it The missing roots $r_4,\,r_3$ and $r_2$ in Proposition 4.26 are,
respectively, the complex dimensions of $G_2$, $A_2$ and $A_1$, namely 14,8 and 3 so that 
$$\eqalign{f_4(s) &= 1/4!\,\,s(s-1)(s-3)(s-14)\cr -f_3(s) &= 1/3!\,\,s(s-1)(s-8)\cr f_2(s) &= 1/2!\,\,
s(s-3)\cr}$$}\vs  Let $k$ be any positive integer. If $m\in \Bbb Z_+$ and $m\geq 2$ let
$C_k(m)$ equal
$C_k$ for the case where $K= SU(m)$. If $m\geq k$ then $k\leq h^{\vee}= m$ and hence, by Theorem
4.23,
$$f_k(m^2-1) = (-1)^k dim\,C_k(m)\eqno (4.70)$$
 In particular note that $f_k(m^2-1)\neq 0$ since $M\geq k$
where $M$ is defined here for $K=SU(m)$. But since $f_k(0) = 0$ and since $f_k$ is a polynomial of
degree $k$ it follows that $f_k(s)$ is determined by the values $f_k(m^2-1)$ for $k$ different
positive values of
$m$. But then (4.70) establishes the following theorem. The result implies that
$f_k(s)$ is encoded in the $k$-dimensional commutative subalgebra structure of $Lie\,Sl(m,\Bbb C)$ for
$k$ different values of
$m$ where
$m\geq k$ and
$m>1$. 
\vs {\bf Theorem 4.28.} {\it Let $k$ be a positive integer. Then $f_k(s)$ is determined by the
numbers $dim\,C_k(m)$ for $k$ different values of $m\in \Bbb Z_+$ where $m\geq k$ and $m>1$.
Furthermore under the assumption $m\in \Bbb Z_+,\,m\geq k$ and $m>1$ one has
$f_k(m^2-1)\neq 0$. }

\vskip 4pt

\centerline{\bf References}\vskip 1.3pc
\parindent=42pt
\rm
\item {[A-F]} R. Adin and A. Frumkin, Rim Hook tableau and Kostant's $\eta$-Function Coefficients,
arXiv:math. CO/0201003 v2
\item {[B]} R. Bott, An Application of the Morse Theory to the Topology of
Lie Groups, {\it Bull. Soc. Math. France}, {\bf 84}(1956), 251-281
\item {[C-P]} P. Cellini and P. Papi, ad-Nilpotent Ideals of a Borel Subalgebra,
{\it J. Alg.}, {\bf 225}(2000), 130-141
\item{[F]} H. D. Fegan, The heat equation and modular forms, {\it J. 
Diff. Geom.}, {\bf 13}(1978), 589-602
\item{[G]} H. Garland, Dedekind's $\eta$-function and the cohomology of infinite dimensional Lie
algebras, {\it PNAS, USA}, {\bf 72}(1975), 2493-2495
\item{[G-L]} H. Garland and J. Lepowsky, Lie algebra homology and the Macdonald-Kac formulas,
{\it Inventiones Math.}, {\bf 34}(1976), 37-76
\item{[Ka]} V. Kac, {\it Infinite dimensional Lie algebras}, Cambridge 1990
\item{[Ko-1]} B.Kostant, Lie algebra cohomology and the generalized Borel-Weil theorem, {\it Ann. of
Math.}, {\bf 74} No. 2. (1961), 329-387
\item{[Ko-2]} B. Kostant, Eigenvalues of a Laplacian and Commutative
Lie subalgebras, {\it Topology}, {\bf 3}(1965), 147-159 
\item{[Ko-3]} B. Kostant, On Macdonald's $\eta$-function Formula, the
Lapalacian and Generalized Exponents, {\it Adv. in Math.} {\bf 20} No.2.
(1976), 179-212
\item{[Ko-4]} B. Kostant, The McKay correspondence, the Coxeter
element and representation theory, {\it Soci\'et\'e Math\'ematique de France, Ast\'erisque, hors
s\'erie}, 1985, 209-255
\item{[Ko-5]} B. Kostant, The set of Abelian ideals of a Borel Subalgebra, Cartan decompositions
and Discrete Series Representations, {\it IMRN}, {\bf 5}(1998), 225-252
\item{[Ko-6]} B. Kostant, On $\wedge\g$ for a Semisimple Lie Algebra $\g$, as an Equivariant Module
over the Symmetric Algebra S(\g), {\it Advanced Studies in Pure Mathematics}, {\bf 26}(2000),
Analysis on Homogeneous Spaces and Representation Theory of Lie Groups, 129-144
\item{[Ku-1]} S. Kumar, Geometry of Schubert cells and cohomology of Kac-moody Lie algebras, {\it
J. Diff. Geom.}, {\bf 20}(1984), 389-431
\item{[Ku-2]} S. Kumar, {\it Kac-Moody Groups, their Flag Varities and
Representation Theory}, PM 204, Birkhauser Boston, 2002
\item{[L]} D. H. Lehmer, Ramanujan's function $\tau(n)$, {\it Duke Math J.}, {\bf 30}(1943)
\item{[Ma-1]} I.G. MacDonald, Affine root systems and Dedekind's $\eta$-function, {\it Invent.
Math.} {\bf 15}(1972), 91-143
\item{[Ma-2]} I. G. Macdonald, {\it Symmetric functions and Hall Polynomials}, Oxford 1995, (2nd
Edition)
\item{[Se]} J-P. Serre, {\it A Course in Arithmetic}, GTM, {\bf 7} Springer-Verlag, 1973
\item{[St]} R. P. Stanley, {\it Enumerative Combinatorics}, {\bf 2}, Cambridge, 1999
\item{[Su]} R. Suter, Abelian ideals in a Borel subalgebra of a complex simple Lie algebra,
arXiv:math.RT/0210463 v1,  30 Oct 2002
\item{[Ze]} S. Zelditch, Macdonald's Identities and the Large N Limit of $YM_2$ on the Cylinder, to
appear, {\it Comm. Math. Physics}

\smallskip
\parindent=30pt
\vskip 8pt
\vbox to 60pt{\hbox{Bertram Kostant}
      \hbox{Dept. of Math.}
      \hbox{MIT}
      \hbox{Cambridge, MA 02139}\vskip 6pt
\noindent E-mail kostant@math.mit.edu}

\end

\end

\end